\documentclass[11pt]{amsart}

\usepackage[english]{babel}
\usepackage[utf8x]{inputenc}
\usepackage[T1]{fontenc}
\usepackage{physics}

\usepackage[a4paper,top=3cm,bottom=2cm,left=3cm,right=3cm,marginparwidth=1.75cm]{geometry}

\usepackage{amsmath}
\usepackage{amssymb}
\usepackage{amsthm}
\usepackage{graphicx}
\usepackage{subcaption}
\usepackage{comment}
\usepackage{enumitem}
\usepackage[numbers]{natbib}
\usepackage[colorinlistoftodos]{todonotes}
\usepackage[colorlinks=true, allcolors=blue]{hyperref}
\usepackage{tikz-cd}

\newtheorem{theorem}{Theorem}[section]
\newtheorem{corollary}{Corollary}[section]

\newtheorem{proposition}{Proposition}[section]
\newtheorem{lemma}{Lemma}[section]
\newtheorem{remark}{Remark}[section]
\newtheorem{example}{Example}[section]
\newtheorem*{theorem*}{Theorem}
\newtheorem{definition}{Definition}[section]

\newcommand*\Span{\text{span}}

\begin{document}

\title{Characterizations of ordered self-adjoint operator spaces}

\author{Travis B. Russell}
\address{Department of Mathematics, Dartmouth College, 
Hanover, NH 03755}
\email{travis.b.russell@dartmouth.edu}
\urladdr{https://sites.google.com/site/travisbrussell} 

\begin{abstract} We describe how self-adjoint ordered operator spaces, also called non-unital operator systems in the literature, can be understood as $*$-vector spaces equipped with a matrix gauge structure. We explain how this perspective has several advantages over other notions of non-unital operator systems in the literature. In particular, the category of matrix gauge $*$-vector spaces includes injective objects and a Webster-Winkler-type duality theorem, both of which we show generally fail with other notions of non-unital operator systems. As applications, we characterize those subspaces of operator systems which are kernels of completely positive maps and define a new operator space structure on the matrix ordered dual of an operator system generalizing the classical notion of a base norm space. \end{abstract}

\maketitle

\section{Introduction}

The study of operator spaces, operator systems, and operator algebras has flourished in the past four decades or so thanks in part to the discovery of abstract characterizations of these objects by Ruan \cite{RUAN1988}, Choi-Effros \cite{CHOIEffros1977}, and Blecher-Ruan-Sinclair \cite{BRS90}, respectively. These characterizations can be viewed as ``intrinsic'' in the sense that they can be expressed in terms of simple relations between the algebraic and metric or order structures which theses spaces naturally possess at each matrix level.

About twenty years ago, an axiomatic characterization for self-adjoint operator spaces with an order structure, also called a ``non-unital operator system'', was proposed by Werner \cite{WERNER2002} (see also Ng \cite{Ng} and Karn \cite{karn_2011} for similar work). Werner's characterization takes as an axiom the existence of sufficiently many positive contractive linear functionals from which the norm and order structure can be recovered. Thus, one could view this characterization as ``extrinsic'' in the sense that it is expressed in terms of the properties of an external object, the dual space. For a detailed history of this progression of results, see Section 4.5 of Blecher's paper \cite{Blecher2007}. Indeed, Blecher notes in the third paragraph of Section 4.5 of \cite{Blecher2007} that the above ``extrinsic'' characterization was the only known condition which guarantees the existence of a completely isometric complete order embedding of an abstractly defined ordered vector space onto a self-adjoint subspace of a $C^*$-algebra without the presence of a unit. Nonetheless, Werner's notion of a non-unital operator system has continued to be useful in the literature, playing a central role in recent work of Connes-van Suijlekom \cite{ConnesSuijlekomII} and Kennedy-Kim-Manor \cite{KennedyManorKim}, for example.

In this paper, we propose two intrinsic abstract characterizations for self-adjoint operator spaces. The first we call a \textit{normal matrix ordered operator space}. This is a $*$-vector space possessing an $L^\infty$-matrix norm and a matrix ordering of closed cones satisfying the condition
\[ a \leq b \leq c \quad \text{implies} \quad \|b\| \leq \max(\|a\|,\|c\|) \]
for all self-adjoint elements $a,b,c$. We prove that these intrinsically characterized objects capture the matrix norm and matrix order structure of a concrete operator space. 
\begin{theorem}
Let $\mathcal{V}$ be a normal matrix ordered operator space. Then there exists a completely isometric complete order embedding $\pi: \mathcal{V} \to B(H)$.
\end{theorem}
\noindent Despite the above theorem, we show that the category whose objects are normal matrix ordered operator spaces and whose morphisms are completely positive completely contractive maps is somewhat pathological, in the sense that some desirable properties that hold in the categories of operator spaces and operator systems fail in the this category. In particular, we show that there are no injective objects in this category, making extensions of completely positive completely contractive maps problematic. The lack of injective objects is a consequence of the non-uniqueness of unitizations for normal matrix ordered operator spaces. We show that a normal matrix ordered operator space can admit completely isometric complete order embeddings of codimension one into non-isomorphic operator systems, i.e. multiple non-isomorphic unitizations. The class of completely positive completely contractive maps on a normal matrix ordered operator space that extend to one unitization differs from the those which extend to the other. These difficulties are remedied by the second perspective, which we describe next.

The second abstract characterization we provide, which is the focus of this work, is by means of \textit{proper matrix gauges}. We define a proper matrix gauge on a $*$-vector space $\mathcal{V}$ to be a sequence of functions $\nu_n: M_n(\mathcal{V}) \to [0,\infty)$ which are gauges, i.e.
\[ \nu_n(tx) = t\nu(x) \text{ for all } t \geq 0 \quad \text{and} \quad \nu_n(x+y) \leq \nu_n(x) + \nu_n(y), \]
satisfy the $L^\infty$-conditions
\[ \nu_k(\alpha^* x \alpha) \leq \|\alpha\|^2 \nu_n(x) \quad \text{and} \quad \nu_{n+k}(x \oplus y) = \max(\nu_n(x), \nu_k(y)) \]
for all $x \in M_n(\mathcal{V})$, $y \in M_k(\mathcal{V})$, and $\alpha \in M_{n,k}$, and are proper, meaning $\nu_n(x)=\nu_n(-x)=0$ implies that $x = 0$.

Matrix gauges arise naturally in $B(H)$ as follows. For each self-adjoint $T \in B(H)$, define
\[ \nu_n(T) = \|T_+\| \]
i.e. the gauge is the norm of the positive part of $T$, or equivalently
\[ \nu_n(T) = \inf \{ t > 0: T \leq tI_H \}. \]
From this gauge, the order and norm structure of an operator space can be uniquely determined, since $T \geq 0$ if and only if $\nu_n(-T) = 0$ and $\|T\| = \max(\nu(T), \nu(-T))$. We prove
\begin{theorem}
Let $\{\nu_n\}$ be a matrix gauge on a $*$-vector space $\mathcal{V}$. Then there exists a Hilbert space $H$, a linear map $\pi: \mathcal{V} \to B(H)$ satisfying $\nu_n(x) = \| \pi(x)_+ \|$ for all self-adjoint $x \in M_n(\mathcal{V})$. Moreover, $\pi$ is a completely isometric complete order embedding with respect to the order given by $x \geq 0$ when $\nu(-x)=0$ and the norm given by $\|x\| = \max(\nu(x),\nu(-x))$.
\end{theorem}
\noindent Therefore matrix gauges yield an alternative abstract characterization for non-unital operator systems.

While normal matrix ordered operator spaces and matrix gauges both abstractly characterize non-unital operator systems, we find that these characterizations are far from equivalent. We consider the category of $*$-vector spaces with matrix gauges, whose morphisms are completely gauge contractive maps, i.e. maps satisfying $\nu(f(x)) \leq \nu(x)$. We show that these maps are always completely positive and completely contractive with respect to the norm and order structure defined by the gauge. Moreover, this category includes $B(H)$ as an injective object for any Hilbert space $H$, in contrast to the category of normal matrix ordered operator spaces. Hence, matrix gauges provide a means for determining precisely which completely positive completely contractive maps admit completely positive completely contractive extensions. We also find that operator systems carry a unique matrix gauge structure, and that every matrix gauge $*$-vector space admits a gauge isometric unitization into a unique operator system. For these reasons, we regard matrix gauges as the more natural way to abstractly characterize non-unital operator systems.

As further justification for this matrix gauge perspective, we provide several non-trivial applications. First, we characterize those subspaces of an operator system which arise as kernels of completely positive maps, clarifying issues raised by Kavruk-Paulsen-Tomforde-Todorov in \cite{Kavruk}. Then we provide a simplication of the Webster-Winkler duality for non-unital operator systems recently developed by Kennedy-Kim-Manor in \cite{KennedyManorKim}. These simplifications are made possible by regarding non-unital operator systems as $*$-vector spaces with matrix gauges rather than normal matrix ordered operator spaces and making use of injective objects. Finally, we use matrix gauges to develop operator representations of duals of operator systems, generalizing the classical notion of a base-norm space. We show, using gauges, that when $\mathcal{V}$ is the matrix ordered Banach space dual of an operator system, there exists a compete order embedding $\pi: \mathcal{V} \to B(H)$ which is completely isometric when restricted to positive elements. Moreover, the matrix norm induced by this representation is equivalent, though not equal, to the operator space dual matrix norm for non-positive elements.

This paper is a significant update of a preprint which first appeared on arXiv in 2015. Most results found in Section 3 and Section 4 follow from the more general paper \cite{RUSSELL2017}, which addresses the case of possibly non-self-adjoint ordered operator spaces. For this reason, the author initially chose not to publish the preprint. However, non-unital operator systems have continued to play an important role in the literature recently (c.f. \cite{ConnesSuijlekomI}, \cite{ConnesSuijlekomII} and \cite{KennedyManorKim}), so the author felt an update of the original manuscript may be of use to the community. The results of Section 5 appear in the original preprint, but have never been published. The results of Section 6 and Section 7 are new. The writing throughout this paper has been updated from the preprint for improved clarity, particularly with regard to the distinction between vector spaces with matrix gauges and Werner's non-unital operator systems. The notation and terminology has also been updated from the preprint in an effort to be more consistent with the literature on operator systems. Finally, we have omitted the preliminary results on non-unital function systems which appeared in the original preprint in order to more efficiently address the noncommutative case.

\section{Preliminaries}

In this section, we recall definitions and fundamental results concerning operator spaces, operator systems, and Werner's non-unital operator system. Throughout, we let $M_n$ denote the $n \times n$ complex matrices, $M_{n,k}$ denote the $n \times k$ rectangular matrices, and $\alpha^* \in M_{k,n}$ denote the conjugate transpose of a matrix $\alpha \in M_{n,k}$.

A \textbf{concrete operator space} is a vector subspace $\mathcal{X}$ of $B(H)$ for some Hilbert space $H$. An \textbf{abstract operator space} is a complex vector space $\mathcal{X}$, together with a sequence $\{\|\cdot\|_n\}_{n=1}^\infty$, satisfying the following properties:
\begin{enumerate}
    \item $\|\cdot\|_n$ is a norm on $\mathcal{X}$ for every $n \in \mathbb{N}$.
    \item $\| x \oplus y \|_{n+m} = \max(\|x\|_n, \|y\|_m)$ whenever $x \in M_n(\mathcal{X})$ and  $y \in M_m(\mathcal{X})$ for every $n,m \in \mathbb{N}$.
    \item $\|\alpha x \beta\|_n \leq \|\alpha\| \|x\|_k \|\beta\|$ whenever $\alpha, \beta^* \in M_{n,k}$ and $x \in M_k(\mathcal{X})$ for all $n,k \in \mathbb{N}$.
\end{enumerate}
In general, a sequence of norms $\{\|\cdot\|_n\}$ satisfying the three conditions above is called an \textbf{$L^\infty$-matrix norm}.

Let $\mathcal{X}$ and $\mathcal{Y}$ be abstract operator spaces, and let $\phi: \mathcal{X} \to \mathcal{Y}$ be a linear map. We let $\phi^{(n)}: M_n(\mathcal{X}) \to M_n(\mathcal{Y})$ denote the inflation map wich applies $\phi$ to the entries of matrices in $M_n(\mathcal{X})$. The map $\phi$ is called \textbf{completely contractive} if $\|\phi^{(n)}(x)\|_n \leq \|x\|_n$ for every $x \in M_n(\mathcal{X})$. If $\phi$ is a scalar multiple of a completely contractive map, then $\phi$ is called \textbf{completely bounded}. If $\|\phi^{(n)}(x)\|_n = \|x\|_n$ for every $x \in M_n(\mathcal{X})$, then $\phi$ is called \textbf{completely isometric} or a \textbf{complete isometry}.

Every concrete operator space is an example of an abstract operator space, since a concrete operator space $\mathcal{X} \subseteq B(H)$ can be equipped with the sequence of norms $\|\cdot\|_n$ given by the operator norm on $M_n(\mathcal{X})$ identified as a subspace of $B(H^{n})$. The following theorem of Ruan shows that every abstract operator space can be realized as a concrete operator space.

\begin{theorem}[Ruan, \cite{RUAN1988}] \label{thm: Ruan}
Let $\mathcal{X}$ be an abstract operator space. Then there exists a Hilbert space $H$ and a complete isometry $\phi: \mathcal{X} \to B(H)$.
\end{theorem}

Ruan's results generalized earlier work of Choi-Effros on operator systems, which we describe now. By a \textbf{$*$-vector space}, we mean a complex vector space $\mathcal{V}$ equipped with a conjugate linear involution $x \mapsto x^*$. The involution extends to $M_n(\mathcal{V})$ by setting $A^* = (a_{ji}^*)_{ij}$ for each $A = (a_{ij}) \in M_n(\mathcal{V})$. We let $M_n(\mathcal{V})_h$ denote the \textbf{hermitian} or \textbf{self-adjoint} elements of $M_n(\mathcal{V})$, i.e. elements satisfying $A = A^*$. By a \textbf{matrix ordered $*$-vector space}, we mean a $*$-vector space $\mathcal{V}$, together with a sequence of cones $\mathcal{C} = \{\mathcal{C}_n\}_{n=1}^\infty$, satisfying the following properties:
\begin{enumerate}
    \item $\mathcal{C}_n$ is a proper cone in $M_n(\mathcal{V})_h$ for every $n \in \mathbb{N}$.
    \item $\alpha^* x \alpha \in \mathcal{C}_k$ whenever $\alpha \in M_{k,n}$ and $x \in \mathcal{C}_k$.
\end{enumerate}
We often use the shorthand $x \geq 0$ to mean $x \in \mathcal{C}_n$ whenever $x \in M_n(\mathcal{V})$ for some matrix ordered $*$-vector space $\mathcal{V}$.

Let $\mathcal{V}$ and $\mathcal{W}$ be a pair of matrix ordered $*$-vector spaces and let $\phi: \mathcal{V} \to \mathcal{W}$ be a linear self-adjoint map. Then $\phi$ is called \textbf{completely positive} if $\phi^{(n)}(x) \geq 0$ whenever $x \geq 0$ and $x \in M_n(\mathcal{V})$. If $\phi$ is injective and both $\phi$ and $\phi^{-1}: \phi(\mathcal{V}) \to \mathcal{V}$ are completely positive, then $\phi$ is called a \textbf{complete order embedding}. A bijective complete order embedding is called a \textbf{complete order isomorphism}.

Let $(\mathcal{S}, \mathcal{C})$ be a matrix ordered $*$-vector space, and suppose that $e \in \mathcal{S}_h$. The triple $(\mathcal{S}, \mathcal{C}, e)$ is called an \textbf{abstract operator system} if $I_n \otimes e$ is an \textbf{archimedean matrix order unit} for the ordered vector space $(M_n(\mathcal{S})_h, \mathcal{C}_n)$; i.e., if for every $x \in M_n(\mathcal{S})_h$ there exists a $t > 0$ such that $x + tI_n \otimes e \geq 0$, and if the relation $x + tI_n \otimes e \geq 0$ for all $t > 0$ implies that $x \geq 0$. 

We remark that every abstract operator system is an example of an abstract operator space. Indeed, if $(\mathcal{S}, \mathcal{C}, e)$ is an abstract operator system, then the formula
\begin{equation} \label{eqn: order unit norm} \|x\|_n^e := \inf \{ t > 0 : \begin{pmatrix} 0 & x \\ x^* & 0 \end{pmatrix} \leq t I_{2n} \otimes e \} \end{equation}
defines an $L^\infty$-matrix norm on the vector space $\mathcal{S}$.

A \textbf{concrete operator system} is a self-adjoint vector subspace of $B(H)$ containing the identity operator $I$. Given a concrete operator system $\mathcal{S} \subseteq B(H)$, we see that $\mathcal{S}$ is an abstract operator system when equipped with the operator adjoint from $B(H)$, the positive cone given by the identification of $M_n(\mathcal{S})$ as a subspace of $B(H^n)$, and the unit given by $I$. Morevoer, the $L^\infty$-matrix norm given in Equation (\ref{eqn: order unit norm}) coincides with the operator norm on $B(H^n)$. The following theorem of Choi-Effros shows that every abstract operator system can be identified with a concrete operator system as well.

\begin{theorem}[Choi-Effros \cite{CHOIEffros1977}] \label{thm: Choi-Effros}
Let $(\mathcal{S},\mathcal{C},e)$ be an abstract operator system. Then there exists a Hilbert space $H$ and a unital complete order embedding $\phi: \mathcal{S} \to B(H)$.
\end{theorem}

In light of Theorem \ref{thm: Ruan} and Theorem \ref{thm: Choi-Effros}, we will sometimes ignore the adjectives ``abstract'' and ``concrete'' when referring to operator spaces and operator systems for the remainder of the paper, clarifying when necessary.

We conclude this section by discussing Werner's non-unital operator system. In \cite{WERNER2002}, Werner defines a \textbf{matrix ordered operator space} to be a triple $(\mathcal{V}, \mathcal{C}, \{\|\cdot\|_n\})$, where $\mathcal{V}$ is a $*$-vector space, $\{\|\cdot\|_n\}_{n=1}^\infty$ is an $L^\infty$ matrix norm on $\mathcal{V}$ such that the adjoint map $x \mapsto x^*$ is a complete isometry, and $\mathcal{C}$ is a matrix ordering on $\mathcal{V}$ such that $\mathcal{C}_n$ is closed for each $n \in \mathbb{N}$. Given a matrix ordered operator space $\mathcal{V}$, Werner defines its \textbf{modified numerical radius} to be the quantity
\[ r_n(x) := \sup \{ \phi \begin{pmatrix} 0 & x \\ x^* & 0 \end{pmatrix} : \phi: M_n(\mathcal{V}) \to \mathbb{C} \text{ is completely positive and completely contractive} \}. \]
Werner proved the following.

\begin{theorem}[Werner, \cite{WERNER2002}] \label{thm: Werer's theorem}
Let $(\mathcal{V}, \mathcal{C}, \{\|\cdot\|_n\}_{n=1}^\infty)$ be a matrix ordered operator space. Then there exists a Hilbert space $H$ and a completely contractive order isomorphism $\phi: \mathcal{V} \to B(H)$. Moreover, if $\|\cdot\|_n = r_n(\cdot)$ for every $n \in \mathbb{N}$, then  there exists a Hilbert space $H$ and a completely isometric complete order isomorphism $\phi: \mathcal{V} \to B(H)$.
\end{theorem}

In light of the above result, Werner defines a \textbf{non-unital operator system} to be a matrix ordered operator space $\mathcal{V}$ for which there exists a $k > 0$ such that $k\|x\|_n \leq r_n(x)$ for for all $x \in M_n(\mathcal{V})$. In Remark 4.14 of \cite{WERNER2002}, Werner says that this condition is equivalent to the condition that whenever $u \leq x \leq v$ we have $\|x\|_n \leq k \max(\|u\|_n, \|v\|_n)$, suggesting that a proof will appear in later work --- though we could not find the proof in a follow-up paper. We will verify this claim in the case $k=1$ in Section \ref{sec: Normal ord op spaces}. When $k\|x\|_n \leq r_n(x)$ for $k=1$, we obtain the equality $\|x\|_n = r_n(x)$, so that Theorem \ref{thm: Werer's theorem}, together with these observations, yeilds a characterization of those matrix ordered operator spaces which can be realized as self-adjoint ordered operator subspaces of $B(H)$. In general, we will call a matrix ordered operator space \textbf{normal} if $\|x\|_n \leq \max(\|u\|_n, \|v\|_n)$ whenever $u \leq x \leq v$.

\begin{remark}
\emph{Werner obtains his results by proving that every non-unital operator system $\mathcal{V}$ admits a \textit{unitization}. This consists of a (unital) operator system $\widehat{\mathcal{V}}$ and a completely contractive complete order embedding $i: \mathcal{V} \to \widehat{\mathcal{V}}$ (see Definition 4.7 of \cite{WERNER2002} for details). The order structure on $\widehat{\mathcal{V}}$ is determined by the sets of positive contractive linear functionals on $M_n(\mathcal{V})$ for each $n \in \mathcal{V}$. The embedding $i$ is completely isometric if and only if the modified numerical radius agrees with the given norm (see Corollary 4.12 of \cite{WERNER2002}). Some situations in which $i$ is completely ismometric are discussed in Proposition 4.16 of \cite{WERNER2002}, but a characterization of all such situations is not provided in that paper.} 
\end{remark}

\begin{remark}
\emph{An example of a non-unital operator system where the map $i$ from the previous remark is not completely isometric (or even just isometric) is provided by Blecher, Kirkpatrick, Neal, and Werner in Proposition 1.1 of \cite{BlecherKirkpartickNealWerner2007}. Take $\mathcal{V}$ to be the operator space dual of a C*-algebra $\mathcal{A}$. A matrix $F=(f_{ij}) \in M_n(\mathcal{V})$ is interpreted as a map $F:\mathcal{A} \to M_n$ by setting the $ij$ entry of $F(x)$ equal to $f_{ij}(x)$. We then obtain a matrix ordering on $\mathcal{V}$ by taking $\mathcal{C}_n$ to be the set of completely positive maps $F: \mathcal{A} \to M_n$ and a matrix norm by setting $\|F\|_n = \sup \{ \|F(x)\| : \|x\| \leq 1 \}$. Blecher observes that the modified numerical radius cannot equal the matrix norm in this case. One reason for this failure is because $\mathcal{V}$ is not \textit{normal}, as defined above. For example, if $f \in \mathcal{V}_h$, then there exist $f_+ , f_- \in \mathcal{C}_1$ such that $f = f_+ - f_-$ and $\|f\| = \|f_+\| + \|f_-\|$. Unless $f$ is equal to either $f_+$ or $-f_-$, it is not possible for $\|f\| \leq \max(\|f_+\|, \|f_-\|)$ even though $-f_- \leq f \leq f_+$.}
\end{remark}

\section{The category of matrix gauge $*$-vector spaces} \label{sec: matrix gauges}

In this section, we introduce matrix gauge $*$-vector spaces and explain how they are realized as matrix ordered self-adjoint subspaces of $B(H)$.

\begin{definition}
Let $\mathcal{V}$ be a $*$-vector space. By a \textbf{matrix gauge} on $\mathcal{V}$, we mean a sequence of functions $\{ \nu_n: M_n(\mathcal{V})_h \to [0,\infty) \}_{n=1}^\infty$ which satisfy the following conditions.
\begin{enumerate}
    \item If $x,y \in M_n(\mathcal{V})_h$, then $\nu_n(x+y) \leq \nu_n(x) + \nu_n(y)$.
    \item If $x \in M_n(\mathcal{V})_h$ and $\alpha \in M_{n,k}$, then $\nu_k(\alpha^* x \alpha) \leq \|\alpha\|^2 \nu_n(x)$.
    \item If $x \in M_n(\mathcal{V})_h$ and $y \in M_m(\mathcal{V})_h$, then $\nu_{n+m}(x \oplus y) = \max(\nu_n(x), \nu_n(y))$.
\end{enumerate}
We say that a matrix gauge is \textbf{proper} if $\nu_n(x)=\nu_n(-x)=0$ implies that $x = 0$. For brevity, we sometimes write $\nu$ in place of $\{\nu_n\}_{n=1}^\infty$. When $\nu$ is a proper matrix gauge on a $*$-vector space $\mathcal{V}$, we call the pair $(\mathcal{V}, \nu)$ a \textbf{matrix gauge $*$-vector space}.
\end{definition}

\begin{remark}
\emph{\begin{enumerate}[label=\alph*)]
    \item When $\nu$ is a matrix gauge on $\mathcal{V}$, the functions $\nu_n$ are not necessarily norms. However, they do satisfy $\nu_n(tx) = t \nu_n(x)$ for all $t \geq 0$, since the equalities $tx = (t^{1/2} I_n) x (t^{1/2}I_n)$ for all $t \geq 0$ and $x = (t^{-1/2} I_n) (tx) (t^{-1/2} I_n)$ for all $t > 0$, together with property (2), imply that $\nu_n(tx) \leq t \nu_n(x)$ for all $t \geq 0$ and $t \nu_n(x) \leq \nu_n(tx)$ when $t > 0$.
    \item Matrix gauges were first introduced by Effros and Winkler in \cite{EFFROSWinkler1997}, though their definition includes the possibility that $\nu_n(x)=\infty$. The notion of a proper matrix gauge is our own.
    \item A proper matrix gauge can be extended to a sequence $\nu_n: M_n(\mathcal{V}) \to [0,\infty)$ by setting $\nu_n(x) = \nu_n(\Re(x))$, and the resulting sequence of functions satisfies the $\mathbb{C}$-proper condition of \cite{RUSSELL2017}, namely that the equalities $\nu_n(\pm x) = \nu_n(\pm ix) = 0$ imply that $x = 0$.
\end{enumerate}}
\end{remark}

Before discussing the utility of matrix gauges, we wish to present a few motivating examples.

\begin{example} \label{example: operator space}
\emph{Let $\mathcal{V}$ be a $*$-vector space and suppose that $\{\|\cdot\|_n\}_{n=1}^\infty$ is an $L^\infty$-matrix norm on $\mathcal{V}$. Then $\nu_n(x) := \|x\|_n$ defines a proper matrix gauge on $\mathcal{V}$.}
\end{example}

\begin{example} \label{example: operator system}
\emph{Let $(\mathcal{S}, \mathcal{C}, e)$ be an operator system. Then
\[ \nu_n^e(x) := \inf \{t > 0: x \leq t I_n \otimes e \} \]
defines a proper matrix gauge on $\mathcal{V}$. Indeed, property (1) follows from the implication $x \leq t I_n \otimes e$ and $y \leq r I_n \otimes e$ implies $x+y \leq (t+r) I_n \otimes e$, property (2) follows from the implication $x \leq \|\alpha\|^2 I_n \otimes e$ implies $\alpha^* x \alpha \leq \alpha^* \alpha \otimes e \leq \|\alpha\|^2 I_n \otimes e$, and property (3) follows from fact that $x \oplus y \leq t I_{n+m} \otimes e$ if and only if $x \leq t I_n \otimes e$ and $y \leq t I_n \otimes e$. That $\nu^e$ is proper follows from the Archimedean property, which implies $x \leq t I_n \otimes e$ for all $t > 0$ if and only if $x \leq 0$. Thus $\nu_n^e(x)=\nu_n^e(-x)=0$ implies that $x,-x \geq 0$. Since $\mathcal{C}$ is proper, $x = 0$ in this case. Furthermore, we can recover the matrix ordering $\mathcal{C}$ on $\mathcal{S}$ from the matrix gauge $\nu^e$. By the archimedean property, $x \in \mathcal{C}_n$ if and only if $-x \leq t I_n \otimes e$ for all $t > 0$. Thus, $x \in \mathcal{C}_n$ if and only if $\nu_n^e(-x)=0$. We can also recover the order unit norm, provided in Equation (\ref{eqn: order unit norm}), since 
\[ \| x \|_n^{e} = \max( \nu_{2n}^e\begin{pmatrix} 0 & x \\ x^* & 0 \end{pmatrix}, \nu_{2n}^e \begin{pmatrix} 0 & -x \\ -x^* & 0 \end{pmatrix} ). \]
When $x=x^*$, we may simply take $\|x\|_n^e = \max(\nu_n^e(x), \nu_n^e(-x))$. To see this, notice that whenever $x=x^*$ and
\begin{equation} \label{eqn: matrix norm example}
    \begin{pmatrix} te & \pm x \\ \pm x & te \end{pmatrix} \geq 0
\end{equation}
we have $te \pm x \geq 0$, since conjugating the expression in Inequality (\ref{eqn: matrix norm example}) by $\begin{pmatrix} \frac{1}{\sqrt{2}} & \frac{1}{\sqrt{2}} \end{pmatrix}$ yields the quantity $te \pm x$. On the other hand, if $x=x^*$ and $te \pm x \geq 0$, then we have
\[ \pm \begin{pmatrix} x & x \\ x & x \end{pmatrix}  + t \begin{pmatrix} e & e \\ e & e \end{pmatrix} \geq 0 \quad \text{and} \quad \mp \begin{pmatrix} x & -x \\ -x & x \end{pmatrix}  + t \begin{pmatrix} e & -e \\ -e & e \end{pmatrix} \geq 0. \]
Averaging these expressions yields Inequality (\ref{eqn: matrix norm example}).}
\end{example}

\begin{example} \label{example: C*algebra}
\emph{Let $\mathcal{A}$ be a C*-algebra. For each $x \in M_n(\mathcal{A})_h$, define $\nu_n(x) = \|x_+\|_n$, where $x_+$ is the positive part of $x$; i.e. $x_+ = \frac{1}{2}(x + |x|)$. Then $\nu_n$ defines a matrix gauge on $\mathcal{A}$. Moreover, we have that $x \geq 0$ if and only if $\nu_n(-x) = 0$ and $\|x\|_n = \max(\nu_n(x),\nu_n(-x))$ whenever $x = x^*$. These claims can be verified by observing that $\nu$ is the restriction of the matrix gauge $\nu^e$ on the C*-algebraic unitization of $\mathcal{A}$.}
\end{example}

\begin{example} \label{example: Dual operator system}
\emph{Let $(\mathcal{S},\mathcal{C},e)$ be an operator system. Let $\mathcal{S}^d$ denote the set of bounded linear functionals on $\mathcal{S}^d$. Then $\mathcal{S}^d$ is a matrix ordered $*$-vector space when equipped with the adjoint defined by $f^*(x) := \overline{f(x^*)}$ and the matrix ordering defined by $(f_{ij}) \geq 0$ if and only if the mapping $(f_{ij}): \mathcal{S} \to M_n$ given by $(f_{ij})(x) := (f_{ij}(x))_{ij}$ is completely positive. See Section 4 of \cite{PaulsenTodorovTomforde2011OSS} for more details regarding the ordered vector space $\mathcal{S}^d$. By a result of Choi and Effros, $\mathcal{S}^d$ is an operator system whenever $\mathcal{S}$ is finite dimensional, though the choice of order unit $e^d \in \mathcal{S}^d$ is not unique. When $\mathcal{S}$ is not finite dimensional, $\mathcal{S}^d$ may not be an operator system. However, we will show that $\mathcal{S}^d$ can be regarded as a matrix gauge $*$-vector space. For each $f \in M_n(\mathcal{S}^d)_h$, define
\[ \nu_n(f) := \inf \{ t > 0 : f \leq t s \text{ for some } s \geq 0 \text{ satisfying } s(I_n \otimes e) = I_n \}. \]
We will show in Section \ref{Sec: duals of operator systems} that $(\mathcal{S}^d, \nu)$ is a matrix gauge $*$-vector space and study its properties.}
\end{example}

In light of Examples (\ref{example: operator system}) and (\ref{example: C*algebra}), we make the following definitions.

\begin{definition}
Let $(\mathcal{V}, \nu)$ be a matrix gauge $*$-vector space. For each $n \in \mathbb{N}$, we define
\[ \mathcal{C}_n^{\nu} := \{ x \in M_n(\mathcal{V})_h : \nu_n(-x) = 0\} \]
and for each $x \in M_n(\mathcal{V})$, we define
\[ \| x \|_n^{\nu} = \max( \nu_{2n}\begin{pmatrix} 0 & x \\ x^* & 0 \end{pmatrix}, \nu_{2n} \begin{pmatrix} 0 & -x \\ -x^* & 0 \end{pmatrix}). \]
We let $\mathcal{C}^{\nu}$ denote the sequence $\{ \mathcal{C}_n^{\nu} \}_{n=1}^\infty$ and we let $\|\cdot\|^{\nu}$ denote the sequence $\{ \|\cdot\|_n^{\nu}\}_{n=1}^\infty$. We call $\mathcal{C}^{\nu}$ (respectively $\|\cdot\|^{\nu}$) the matrix ordering (respectively matrix norm) \textbf{induced} by $\nu$.
\end{definition}

\begin{remark}
Whenever $x=x^*$, we have $\|x\|_n^{\nu} = \max(\nu_n(x),\nu_n(-x))$. This is because of the inequalities
\[ \nu_n(x) = \nu_n( \begin{pmatrix} \frac{1}{\sqrt{2}} & \frac{1}{\sqrt{2}} \end{pmatrix} \begin{pmatrix} 0 & x \\ x & 0 \end{pmatrix} \begin{pmatrix} \frac{1}{\sqrt{2}} \\ \frac{1}{\sqrt{2}} \end{pmatrix} ) \leq \nu_{2n} \begin{pmatrix} 0 & x \\ x & 0 \end{pmatrix}, \]
\[ \nu_n(-x) = \nu_n( \begin{pmatrix} \frac{1}{\sqrt{2}} & \frac{1}{\sqrt{2}} \end{pmatrix} \begin{pmatrix} 0 & -x \\ -x & 0 \end{pmatrix} \begin{pmatrix} \frac{1}{\sqrt{2}} \\ \frac{1}{\sqrt{2}} \end{pmatrix} ) \leq \nu_{2n} \begin{pmatrix} 0 & -x \\ -x & 0 \end{pmatrix}, \]
\[ \nu_{2n} \begin{pmatrix} 0 & x \\ x & 0 \end{pmatrix} = \nu_{2n} \begin{pmatrix} \frac{1}{\sqrt{2}} & \frac{1}{\sqrt{2}} \\ \frac{1}{\sqrt{2}} & -\frac{1}{\sqrt{2}} \end{pmatrix} \begin{pmatrix} x & 0 \\ 0 & -x \end{pmatrix} \begin{pmatrix} \frac{1}{\sqrt{2}} & \frac{1}{\sqrt{2}} \\ \frac{1}{\sqrt{2}} & -\frac{1}{\sqrt{2}} \end{pmatrix} \leq \nu_{2n} \begin{pmatrix} x & 0 \\ 0 & -x \end{pmatrix} = \max(\nu_n(x),\nu_n(-x)) \]
and
\[ \nu_{2n} \begin{pmatrix} 0 & -x \\ -x & 0 \end{pmatrix} = \nu_{2n} \begin{pmatrix} \frac{1}{\sqrt{2}} & \frac{1}{\sqrt{2}} \\ \frac{1}{\sqrt{2}} & -\frac{1}{\sqrt{2}} \end{pmatrix} \begin{pmatrix} -x & 0 \\ 0 & x \end{pmatrix} \begin{pmatrix} \frac{1}{\sqrt{2}} & \frac{1}{\sqrt{2}} \\ \frac{1}{\sqrt{2}} & -\frac{1}{\sqrt{2}} \end{pmatrix} \leq \nu_{2n} \begin{pmatrix} -x & 0 \\ 0 & x \end{pmatrix} = \max(\nu_n(x),\nu_n(-x)). \]
\end{remark}

We will show below that whenever $(\mathcal{V}, \nu)$ is a matrix gauge $*$-vector space, then $(\mathcal{V}, \mathcal{C}^{\nu}, \|\cdot\|^{\nu})$ is a normal matrix ordered operator space and hence admits a completely isometric complete order embedding into $B(H)$ (see Corollary \ref{cor: normal OVS}). Rather than proving all this directly, we will show a stronger result: that every matrix gauge $*$-vector space admits a completely gauge isometric embedding into $B(H)$. To explain this, we must describe completely gauge contractive maps, which we consider to be the natural morphisms in the category of matrix gauge $*$-vector spaces.

\begin{definition}
Let $(\mathcal{V}, \nu)$ and $(\mathcal{W}, \omega)$ be a pair of matrix gauge $*$-vector spaces. Let $\phi: \mathcal{V} \to \mathcal{W}$ be a linear map. Then $\phi$ is called \textbf{completely gauge contractive} if $\omega_n(\phi^{(n)}(x)) \leq \nu_n(x)$ for every $x \in M_n(\mathcal{V})_h$. The map $\phi$ is called \textbf{completely gauge isometric} or a \textbf{complete gauge isometry} if $\omega_n(\phi^{(n)}(x)) = \nu_n(x)$ for every $x \in M_n(\mathcal{V})_h$.
\end{definition}

\begin{proposition} \label{prop: gauge contractive is cp and cc}
Let $(\mathcal{V}, \nu)$ and $(\mathcal{W}, \omega)$ be a pair of matrix gauge $*$-vector spaces. Let $\phi: \mathcal{V} \to \mathcal{W}$ be a self-adjoint linear map. If $\phi$ is completely gauge contractive, then $\phi$ is completely positive and completely contractive as a map from $(\mathcal{V}, \mathcal{C}^{\nu}, \|\cdot\|^{\nu})$ to $(\mathcal{W}, \mathcal{C}^{\omega}, \|\cdot\|^{\omega})$. If $\phi$ is completely gauge isometric, then $\phi$ is a complete order embedding and completely isometric as a map from $(\mathcal{V}, \mathcal{C}^{\nu}, \|\cdot\|^{\nu})$ to $(\mathcal{W}, \mathcal{C}^{\omega}, \|\cdot\|^{\omega})$.
\end{proposition}

\begin{proof}
First, let us assume that $\phi$ is completely gauge contractive. Suppose that $x \in M_n(\mathcal{V})$ and $x \geq 0$. Then since
\[ 0 \leq \omega(-\phi^{(n)}(x)) = \omega(\phi^{(n)}(-x)) \leq \nu(-x) = 0 \]
we have that $\phi^{(n)}(x) \geq 0$. So $\phi$ is completely positive. Furthermore, for any $x \in M_n(\mathcal{V})$,
\begin{eqnarray} 
\| \phi^{(n)}(x) \|_n^{\omega} & = & \max( \omega_{2n} \begin{pmatrix} 0 & \phi^{(n)}(x) \\ \phi^{(n)}(x)^* & 0 \end{pmatrix}, \omega_{2n} \begin{pmatrix} 0 & -\phi^{(n)}(x) \\ -\phi^{(n)}(x)^* & 0 \end{pmatrix} ) \nonumber \\
& = & \max( \omega_{2n} (\phi^{(2n)} \begin{pmatrix} 0 & x \\ x^* & 0 \end{pmatrix}), \omega_{2n} (\phi^{(2n)} \begin{pmatrix} 0 & -x \\ -x^* & 0 \end{pmatrix} )) \nonumber \\
& \leq & \max( \nu_{2n} \begin{pmatrix} 0 & x \\ x^* & 0 \end{pmatrix}, \nu_{2n} \begin{pmatrix} 0 & -x \\ -x^* & 0 \end{pmatrix} ) \nonumber \\
& = & \|x\|_n^{\nu}. \nonumber
\end{eqnarray}
Therefore $\phi$ is completely contractive.

Next, assume that $\phi$ is completely gauge isometric. It follows that $\phi$ is completely gauge contractive, and hence completely positive and completely contractive. We claim that $\phi$ is one-to-one. Indeed, suppose that $\phi(x) = 0$. Then
\begin{eqnarray} 
0 & = & \| \phi^{(n)}(x) \|_n^{\omega} \nonumber \\
& = & \max( \omega_{2n} \begin{pmatrix} 0 & \phi^{(n)}(x) \\ \phi^{(n)}(x)^* & 0 \end{pmatrix}, \omega_{2n} \begin{pmatrix} 0 & -\phi^{(n)}(x) \\ -\phi^{(n)}(x)^* & 0 \end{pmatrix} ) \nonumber \\
& = & \max( \omega_{2n} (\phi^{(2n)} \begin{pmatrix} 0 & x \\ x^* & 0 \end{pmatrix}), \omega_{2n} (\phi^{(2n)} \begin{pmatrix} 0 & -x \\ -x^* & 0 \end{pmatrix} )) \nonumber \\
& = & \max( \nu_{2n} \begin{pmatrix} 0 & x \\ x^* & 0 \end{pmatrix}, \nu_{2n} \begin{pmatrix} 0 & -x \\ -x^* & 0 \end{pmatrix} ). \nonumber
\end{eqnarray}
Since $\nu$ is proper, we must conclude that
\[ \begin{pmatrix} 0 & x \\ x^* & 0 \end{pmatrix} = \begin{pmatrix} 0 & 0 \\ 0 & 0 \end{pmatrix} \]
and hence $x = 0$. Since $\phi$ is one-to-one and completely gauge isometric, the map $\phi^{-1}: \phi(\mathcal{V}) \to \mathcal{V}$ is well-defined and completely gauge isometric. Since $\phi^{-1}$ is completely gauge isometric, it is completely positive and completely contractive. It follows that $\phi$ is a completely isometric complete order embedding.
\end{proof}

Proposition \ref{prop: gauge contractive is cp and cc} states that every completely gauge contractive map is completely positive and completely contractive, but we will see later that the converse of that statement fails generally. In Section \ref{sec: Normal ord op spaces}, we will characterize precisely when the converse holds. For now, we mention just one important case where the converse holds.

\begin{proposition} \label{prop: cpcc maps are cgc on operator systems}
Let $\mathcal{S}$ and $\mathcal{T}$ be operator systems with units $e$ and $f$, respectively, and associated matrix gauges $\nu^e$ and $\nu^f$, respectively (see Example \ref{example: operator system}). Then every completely postive completely contractive map $\phi: \mathcal{S} \to \mathcal{T}$ is completely gauge contractive. Consequently, $\phi$ is completely positive and completely contractive if and only if $\phi$ is completely gauge contractive, and $\phi$ is a complete gauge isometry if and only if $\phi$ is a completely isometric complete order embedding.
\end{proposition}

\begin{proof}
Suppose that $\phi: \mathcal{S} \to \mathcal{T}$ is completely positive and completely contractive, and that $x \in M_n(\mathcal{S})_h$ and $\nu_n^e(x) \leq 1$. Then $x \leq I_n \otimes e$. Thus, $\phi^{(n)}(x) \leq \phi^{(n)}(I_n \otimes e) \leq I_n \otimes f$. It follows that $\nu_n^f(\phi^{(n)}(x)) \leq 1$. We conclude that $\phi$ is completely gauge contractive. The remaining statements follow from Proposition \ref{prop: gauge contractive is cp and cc} and from considering the case when both $\phi$ and $\phi^{-1}$ are completely positive and completely contractive.
\end{proof}

Although we have not yet shown that $\mathcal{C}^{\nu}$ is a matrix ordering or that $\|\cdot\|^{\nu}$ is an $L^\infty$-matrix norm, we have seen in Example \ref{example: operator system} a case where these structures coincide with the matrix norm and matrix ordering of an operator system. In light of Proposition \ref{prop: gauge contractive is cp and cc}, the existence of a gauge isometric mapping from $(\mathcal{V},\nu)$ into an operator system is sufficient to prove that $\mathcal{C}^{\nu}$ is a matrix ordering and that $\|\cdot\|^{\nu}$ is an $L^\infty$-matrix norm. We will show that this can always be done. To this end, we define the unitization of a matrix gauge $*$-vector space.

In the following, we identify the vector spaces $M_n(\mathcal{V} \oplus \mathbb{C})$ and $M_n(\mathcal{V}) \oplus M_n$ in the obvious way. When $X \in (M_n)_{h}$, we let $X \gg 0$ indicate that $X \geq 0$ and $X$ is invertible, and we let $X_t := tI_n - X$ whenever $t > 0$. To simplify notation, we let $\widehat{\mathcal{V}}$ denote the vector space $\mathcal{V} \oplus \mathbb{C}$.

\begin{definition}
Let $(\mathcal{V}, \nu)$ be a matrix ordered $*$-vector space. For each $n \in \mathbb{N}$, we define a map $\widehat{\nu}_n: M_n(\widehat{\mathcal{V}})_h \to [0,\infty)$ by setting
\[ \widehat{\nu}_n(A, X) = \inf \{ t > 0 : X_t \gg 0 \text{ and } \nu_n((X_t)^{-1/2}A(X_t)^{-1/2}) \leq 1 \}. \]
We let $\widehat{\nu}$ denote the sequence $\{\widehat{\nu}_n\}_{n=1}^\infty$. We call the pair $(\widehat{\mathcal{V}}, \widehat{\nu})$ the \textbf{unitization} of $(\mathcal{V},\nu)$.
\end{definition}

We now show that $(\widehat{\mathcal{V}}, \widehat{\nu})$ is an operator system.

\begin{theorem} \label{thm: Unitization}
Let $(\mathcal{V}, \nu)$ be a matrix gauge $*$-vector space and let $(\widehat{\mathcal{V}}, \widehat{\nu})$ be its unitization, and let $e=(0,1) \in \widehat{\mathcal{V}}$. Then $(\widehat{\mathcal{V}}, \mathcal{C}^{\widehat{\nu}}, e)$ is an operator system satisfying $\widehat{\nu} = \nu^e$, i.e.
\[ \widehat{\nu}_n(x) = \inf \{t > 0 : x \leq tI_n \otimes e \} \]
for all $x \in M_n(\widehat{\mathcal{V}})_h$. Moreover, the mapping $i: \mathcal{V} \to \widehat{\mathcal{V}}$ given by $x \mapsto (x,0)$ is completely gauge isometric.
\end{theorem}

\begin{proof}
To see that $i$ is completely gauge isometric, notice that
\[ \widehat{\nu}_n(A,0) = \inf\{ t > 0 : \nu_n(A) \leq t \} \]
since $0_t := tI_n$, and hence $\widehat{\nu}_n(A,0)=\nu_n(A)$ for each $A \in M_n(\mathcal{V})_h$. To show that $\mathcal{C}^{\widehat{\nu}}$ is a matrix cone, it suffices to show that it is closed under conjugation by scalar matrices and direct sums. Suppose that $(A,X) \geq 0$ and $(B,Y) \geq 0$, where $(A,X) \in M_n(\widehat{\mathcal{V}})$ and $(B,Y) \in M_m(\widehat{\mathcal{V}})$. Then $\widehat{\nu}_n(-A,-X) = \widehat{\nu}_m(-B,-Y)=0$. Since \[ (-X)_t^{-1/2}(-A)(-X)_t^{-1/2} \oplus (-Y)_t^{-1/2}(-B)(Y)_t^{-1/2}) = (-X \oplus -Y)_t^{-1/2} (-A \oplus -B) (-X \oplus -Y)_t^{-1/2} \]
and since $\nu_{n+m}((-X)_t^{-1/2}(-A)(-X)_t^{-1/2} \oplus (-Y)_t^{-1/2}(-B)(Y)_t^{-1/2})$ is equal to
\[ \max( \nu_n((-X)_t^{-1/2}(-A)(-X)_t^{-1/2}), \nu_m((-Y)_t^{-1/2}(-B)(Y)_t^{-1/2})) \leq 1 \]
for all $t > 0$, we see that $\widehat{\nu}_{n+m}((-A,-X) \oplus (-B,-Y)) = 0$ and hence $(A,X) \oplus (B,Y) \geq 0$. To see that $\mathcal{C}^{\widehat{\nu}}$ is closed under conjugation, let $\alpha \in M_{n,k}$ and consider $\alpha^*(A,X)\alpha = (\alpha^*A\alpha, \alpha^*X\alpha) \in M_k(\widehat{\mathcal{V}})$. Since $(-X)_t \gg 0$ for all $t > 0$, $(-\alpha^*X\alpha)_t \gg 0$ for all $t > 0$. Observe that
\[ (-\alpha^* X \alpha)_t^{-1/2} \alpha^* (-A) \alpha (-\alpha^* X \alpha)_t^{-1/2} \]
equals \[ (-\alpha^* X \alpha)_t^{-1/2} \alpha^* (-X)_t^{1/2} (-X)_t^{-1/2}(-A)(-X)_t^{-1/2} (-X)_t^{1/2} \alpha (-\alpha^* X \alpha)_t^{-1/2} \]
and hence
\begin{eqnarray} \nu_k((-\alpha^* X \alpha)_t^{-1/2} \alpha^* (-A) \alpha (-\alpha^* X \alpha)_t^{-1/2}) & \leq & \|(-\alpha^* X \alpha)_t^{-1/2} \alpha^* (-X)_t^{1/2}\|^2 \nu_n((-X)_t^{-1/2}(-A)(-X)_t^{-1/2}) \nonumber \\
& \leq & \|(-\alpha^* X \alpha)_t^{-1/2} \alpha^* (-X)_t^{1/2}\|^2. \nonumber
\end{eqnarray}
If $\alpha^* \alpha \leq I_k$, then
\[ (-\alpha^* X \alpha)_t^{-1/2} \alpha^* (-X)_t^{1/2} \alpha (-\alpha^* X \alpha)_t^{-1/2} \leq (-\alpha^* X \alpha)_t^{-1/2} (-\alpha^* X \alpha)_t (-\alpha^* X \alpha)_t^{-1/2} = I_k. \]
Hence
\[ \nu_k((-\alpha^* X \alpha)_t^{-1/2} \alpha^* (-A) \alpha (-\alpha^* X \alpha)) \leq 1 \]
for all $t > 0$.  So $\nu_k(\alpha^*(A,X)\alpha) = 0$ and hence $\alpha^*(A,X)\alpha \geq 0$ in this case. For the general case, replace $\alpha$ with $\alpha/\|\alpha\|$, $X$ with $\|\alpha\|^2 X$ and $A$ with $\|\alpha\|^2 A$.

It remains to check that $e$ is an archimedean order unit and that \[ \widehat{\nu}_n(x) = \inf \{t > 0 : x \leq t I_n \}. \]
Suppose that $(A,X) \in M_n(\widehat{\mathcal{V}})_h$. Then since $X=X^*$, there exists $t > 0$ such that $X_t \gg 0$. Now
\[ \nu_n((X)_t^{-1/2} A (X)_t^{-1/2}) \leq \|(X)_t^{-1/2}\|^2 \nu_n(A) = \|X_t^{-1}\| \nu_n(A). \]
To show that $e=(0,1)$ is an order unit, we must show that there exists $r > 0$ such that $(A,X + rI_n) \geq 0$; i.e. $\widehat{\nu}_n(-A, -X-rI_n) = 0$. Choose $r$ such that $X + rI_n \gg 0$ and such that the smallest eigenvalue of $X + rI_n$, say $\lambda_r$, satisfies $\lambda_r^{-1} \leq \nu_n(-A)^{-1}$ (if $\nu_n(-A) = 0$, there will be nothing to show). Then for every $t > 0$,
\[ \nu_n((-X-rI_n)^{-1/2}(-A)(-X-rI_n)^{-1/2}) \leq \|(-X-rI_n)_t^{-1}\|\nu_n(-A) \leq 1 \]
for all $t > 0$, since $(-X-rI_n)_t = (t+r)I_n + X$. It follows that $e=(0,1)$ is an order unit. We conclude from checking
\[ \widehat{\nu}_n(A,X) = \inf \{t > 0 : (A,X) \leq tI_n \otimes e\} \]
from wence the archimedean property will follow. Suppose that $t > 0$ and that $(0,tI_n) \geq (A,X)$; i.e. $\widehat{\nu}_n(A, X-tI_n) = 0$. Then for every $r > 0$, we have $(X-tI_n)_r \gg 0$ and 
\[ \nu_n((X-tI_n)_r^{-1/2}A(X-tI_n)_r^{-1/2}) \leq 1. \] 
Notice that $(X-tI_n)_r = (r+t)I_n - X = (X)_{t+r}$. It follows that $t \geq u_n(A,X)$. A similar calculation shows that whenever $t > \widehat{\nu}_n(A,X)$, we have $(0,tI_n) \geq (A,X)$. So $\widehat{\nu}_n$ has the desired property. Consequently, we have $\widehat{\nu}_n(-A,-X) = 0$ if and only if $0 \leq (A,X) + tI_n \otimes e$ for every $t > 0$. So $e$ is an Archimedean matrix order unit.
\end{proof}

\begin{corollary} \label{cor: Representations of matrix gauge spaces}
Let $(\mathcal{V}, \nu)$ be a matrix gauge $*$-vector space. Then there exists a Hilbert space $H$ and a map $\phi:\mathcal{V} \to B(H)$ such that \[ \nu_n(x) = \| \phi^{(n)}(x)_+ \| \]
for all $x \in M_n(\mathcal{V})_h$. Moreover, $\phi$ is a completely isometric complete order embedding on the normal matrix ordered operator space $(\mathcal{V}, \mathcal{C}^{\nu}, \|\cdot\|^{\nu})$.
\end{corollary}

\begin{proof}
By Theorem \ref{thm: Unitization}, there exists an operator system $(\mathcal{S}, \mathcal{C}, e)$ and a map $i: \mathcal{V} \to \mathcal{S}$ which is completely gauge isometric when $\mathcal{S}$ is equipped with the matrix gauge $\nu^e$. Since $\mathcal{C} = \mathcal{C}^{\nu^e}$ and $\|\cdot\|^{e} = \|\cdot\|^{\nu^e}$, and since $i$ is completely gauge isometric, $i$ is a completely isometric complete order embedding on $(\mathcal{V}, \mathcal{C}^{\nu}, \|\cdot\|^{\nu})$. By the Theorem \ref{thm: Choi-Effros}, there exists a unital complete order embedding $\pi: \mathcal{S} \to B(H)$. Finally, for all $x \in M_n(\mathcal{V})_h$,
\[ \nu_n(x) = \nu^e_n(i^{(n)}(x)) = \| \pi \circ i^{(n)}(x)_+\|. \]
It follows that $(\mathcal{V}, \mathcal{C}^{\nu}, \|\cdot\|^{\nu})$ is a normal matrix ordered $*$-operator space, since it is completely isometrically order isomorphic to a concrete self-adjoint operator space.
\end{proof}

We will conclude this section with an analogue of Arveson's Extension Theorem for matrix gauge $*$-vector spaces. We first show that completely gauge contractive maps extend to unital completely positive maps on their unitizations.

\begin{proposition} \label{prop: unitize maps}
Let $(\mathcal{V},\nu)$ and $(\mathcal{W},\omega)$ be a pair of matrix gauge $*$-vector spaces, and suppose that $\phi: \mathcal{V} \to \mathcal{W}$ is a self-adjoint linear map. Then the unique unital extension of $\widehat{\phi}: \widehat{\mathcal{V}} \to \widehat{\mathcal{W}}$ given by $\widehat{\phi}(x,\lambda) = (\phi(x),\lambda)$ is completely positive. Moreover, when $\phi$ is completely gauge isometric, then $\widehat{\phi}$ is a unital complete order embedding.
\end{proposition}

\begin{proof}
First, assume that $\phi$ is completely gauge contractive. By Proposition \ref{prop: gauge contractive is cp and cc}, it suffices to check that $\widehat{\phi}$ is completely gauge contractive. To this end, suppose that $(A,X) \in M_n(\widehat{\mathcal{V}})$ with $A=A^*$ and $X=X^*$. Assume that $t > 0$, $X_t \gg 0$, and $\nu_n((X_t)^{-1/2}(A)(X_t)^{-1/2}) \leq 1$. Then
\[ \omega_n((X_t)^{-1/2}(\phi^{(n)}(A))(X_t)^{-1/2}) = \omega_n(\phi^{(n)}((X_t)^{-1/2}(A)(X_t)^{-1/2})) \leq \nu_n((X_t)^{-1/2}(A)(X_t)^{-1/2}) \leq 1. \]
It follows that $\widehat{\omega}_n(\widehat{\phi}^{(n)}(A,X)) \leq t$. Since we may choose a decreasing sequence $t_n$ of values satisfying the above conditions with $t_n \to \widehat{\nu}_n(A,X)$, we conclude that $\widehat{\omega}_n(\widehat{\phi}^{(n)}(A,X)) \leq \widehat{\nu}_n(A,X)$. So $\widehat{\phi}$ is completely gauge contractive by Proposition \ref{prop: cpcc maps are cgc on operator systems}.

Now assume $\phi$ is completely gauge isometric. Then both $\phi$ and the restriction of $\phi^{-1}$ to the range of $\phi$ are completely gauge contractive. By the above, both $\widehat{\phi}$ and $\widehat{\phi}^{-1}$ are completely gauge contractive. It follows that $\widehat{\phi}$ is a complete order embedding in this case. \end{proof}

Proposition \ref{prop: unitize maps} allows us to regard unitization as a functor from the category of matrix gauge $*$-vector spaces, whose morphisms are completely gauge contractive maps, to the category of operator systems, whose morphisms are completely positive maps. This functor sends a matrix gauge $*$-vector space to its unitization and completely gauge contractive maps to their unique unital extensions. The next proposition concerns the behavior of the unitization functor on objects which are already unital. We will say more about this case at the end of Section \ref{sec: kernels and quotients}, where the next result will be used to characterize matrix gauge $*$-vector spaces which are completely gauge isometric to operator systems.

\begin{proposition} \label{prop: unitize operator system}
Let $(\mathcal{V}, \mathcal{C}, e)$ be an operator system and let $\nu^e$ denote the associated matrix gauge. Then the unital map $\pi: \widehat{\mathcal{V}} \to \mathcal{V}$ given by $\pi: (x,\lambda) \mapsto x + \lambda e$
is a unital completely positive map.
\end{proposition}

\begin{proof}
By Proposition \ref{prop: gauge contractive is cp and cc}, it suffices to check that $\pi$ is completely gauge contractive. To this end, suppose that $(A,X) \in M_n(\widehat{\mathcal{V}})$ with $A=A^*$ and $X=X^*$. Assume that $t > 0$, $X_t \gg 0$, and $\nu_n((X_t)^{-1/2}(A)(X_t)^{-1/2}) \leq 1$. Then $(X_t)^{-1/2}(A)(X_t)^{-1/2} \leq I_n \otimes e,$ and hence $A \leq X_t = tI_n - X.$
It follows that $\pi^{(n)}(A,X) = A+X \leq t I_n$ and hence $\nu_n(\pi^{(n)}(A,X)) \leq t$. Since we may choose a decreasing sequence $t_n$ of values satisfying the above conditions with $t_n \to \widehat{\nu}_n(A,X)$, we conclude that $\pi$ is completely gauge contractive.
\end{proof}

\begin{corollary} \label{cor: Unique unital extension op sys}
Let $(\mathcal{V}, \nu)$ be a matrix gauge $*$-vector space and let $\mathcal{S}$ be an operator system. Let $\varphi: \mathcal{V} \to \mathcal{S}$ be completely gauge contractive. Then the unique unital extension $\widetilde{\varphi}: \widehat{\mathcal{V}} \to \mathcal{S}$ is completely positive.
\end{corollary}

\begin{proof}
Proposition \ref{prop: unitize maps} gives a unique unital completely positive map $\widehat{\varphi}: \widehat{\mathcal{V}} \to \widehat{\mathcal{S}}$, and Proposition \ref{prop: unitize operator system} gives a unique unital completely positive $\pi: \widehat{\mathcal{S}} \to \mathcal{S}$. Composing these maps, we obtain the completely positive map $\widetilde = \pi \circ \widehat{\phi}$.
\end{proof}

We now apply these results to prove a non-unital analogue of the Arveson extension theorem. In the following, we regard $B(H)$ as an operator system with its associated matrix gauge, or equivalently with the matrix gauge $\nu_n(T) = \|T_+\|$ (see Example \ref{example: C*algebra}).

\begin{theorem} \label{thm: gauge Arveson}
Let $\mathcal{V} \subseteq {\mathcal{W}}$ be an inclusion of matrix gauge $*$-vector spaces, and let $H$ be a Hilbert space. If $\phi: \mathcal{V} \to B(H)$ is completely gauge contractive, then there exists a completely gauge contractive extension $\widetilde{\phi}: \mathcal{W} \to B(H)$.
\end{theorem}

\begin{proof}
Let $\phi: \mathcal{V} \to B(H)$ be completely gauge contractive. By Proposition \ref{prop: unitize maps}, the unital extension of the inclusion map $\widehat{\mathcal{V}} \subseteq \widehat{\mathcal{W}}$ is a unital complete order embedding of operator systems. Let $\widehat{\phi}: \widehat{\mathcal{V}} \to \widehat{B(H)}$ be the unique unital completely positive extension of $\phi$ from the unitization of $\mathcal{V}$ to the unitization of $B(H)$. Let $\pi: \widehat{B(H)} \to B(H)$ be the unital completely positive map described in Proposition \ref{prop: unitize operator system}. Then $\pi \circ \widehat{\phi}: \widehat{\mathcal{V}} \to B(H)$ is a unital completely positive extension of $\phi$ to $\widehat{\mathcal{V}}$. By Arveson's extension theorem, there exists a completely positive extension $\widetilde{\phi}: \widehat{\mathcal{W}} \to B(H)$. Since $\widetilde{\phi}$ is completely positive and unital, it is completely gauge contractive by Proposition \ref{prop: cpcc maps are cgc on operator systems}. The restriction of this extension to $\mathcal{W}$ is a completely gauge contractive extension of $\phi$.
\end{proof}

\begin{figure}[h!]
\[
\begin{tikzcd}
\widehat{\mathcal{W}} &  \arrow[hook, l] \mathcal W \arrow[rd, "\widetilde{\phi}" black] & \widehat{B(H)} \arrow[d, "\pi"] \\
\arrow[dashrightarrow, urr, "\widehat{\phi}"] \arrow[hook, u] \widehat{\mathcal{V}}  & \arrow[hook,l] \mathcal V \arrow[r, "\phi" black] & B(H)
\end{tikzcd}
\]
\caption{Map of the proof of Theorem \ref{thm: gauge Arveson}}
\end{figure}

We should emphasize again that while Proposition \ref{prop: gauge contractive is cp and cc} implies that completely gauge contractive maps are always completely positive and completely contractive, we will see in the next section that the converse can fail. Thus Theorem \ref{thm: gauge Arveson} does not imply that every completely positive completely contractive map on a self-adjoint operator space has a completely positive completely contractive extension. The next section addresses this issue in detail, starting with some examples.

\section{The category of normal matrix ordered operator spaces} \label{sec: Normal ord op spaces}

In this section, we will explain some subtle differences between matrix gauge $*$-vector spaces and Werner's non-unital operator spaces. This discussion will lead naturally to the question of when a completely positive completely contractive map has a completely positive completely contractive extension. We begin with an illustrative example.

\begin{example} \label{example: non-gauge isometric spaces}
\emph{Let $\mathcal{V}$ denote the one-dimensional self-adjoint operator space spanned by the matrix
\[ A = \begin{pmatrix} 1 & 0 \\ 0 & -1/2 \end{pmatrix} \]
and let $\mathcal{W}$ be the one-dimensional self-adjoint operator space spanned by the matrix
\[ B = \begin{pmatrix} 1 & 0 \\ 0 & -1 \end{pmatrix}. \]
Then $\mathcal{V}$ and $\mathcal{W}$ inherit the same operator norm and matrix ordering from $M_2$, namely
\[ \|X \otimes A\|_n = \|X \otimes B\| = \|X\| \]
for all $X \in M_n$, and $\mathcal{C}_n = \{0_n\}$ for each $n \in \mathbb{N}$. However, the restriction $\nu$ of $\nu^e$ on $M_2$ to $\mathcal{V}$ gives $\nu_n(A) = 1$ and $\nu_n(-A)=1/2$, while the restriction $\omega$ of $\nu^e$ on $M_2$ to $\mathcal{W}$ gives $\omega_n(B) = \omega_n(-B)=1$. Thus the mapping $\lambda A \mapsto \lambda B$ from $\mathcal{V}$ to $\mathcal{W}$ is a completely isometric complete order isomorphism, but it is not gauge contractive.}

\emph{Let $\mathcal{V}'$ and $\mathcal{W}'$ denote the unitizations of the operator spaces $\mathcal{V}$ and $\mathcal{W}$ in $M_2$, respectively, i.e. $\mathcal{V}' = \Span\{A,I_2\}$ and $\mathcal{W}' = \Span\{B, I_2\}$. Define $\varphi: \mathcal{V} \to \mathbb{C}$ by $\varphi(A)=-1$ and $\psi: \mathcal{W} \to \mathbb{C}$ by $\psi(B)=-1$. Then both $\varphi$ and $\psi$ are easily seen to be completely positive completely contractive maps. Moreover, $\psi$ admits a ucp extension to $\mathcal{W}'$ given by $\psi'(\lambda B + \mu I_2) = \mu-\lambda$. In fact, this is the vector state corresponding to the canonical basis vector $(0,1)^T = \vec{e}_2 \in \mathbb{C}^2$. However, $\varphi$ does not extend to a ucp map on $\mathcal{V}'$. This is because the numerical range of the matrix $A$ is $[-1/2,1]$, which does not contain $\varphi(A)$. Thus, no extension of $\varphi$ can be ucp, since states map operators to values within their numerical range. In fact, the only unital extension of $\varphi$ is given by $\varphi'(\lambda A + \mu I_2) = \mu - \lambda$, but $\varphi'$ has the property}
\[ \varphi' \begin{pmatrix} 3/2 & 0 \\ 0 & 0 \end{pmatrix} = (1/2)-1 = -1/2. \]
\emph{Thus, $\varphi'$ is not a positive functional on $\mathcal{V}'$.}
\end{example}

We will explain later some of the consequences of the above example. For now, we note that Example \ref{example: non-gauge isometric spaces} shows that a given matrix ordered $*$-operator space $(\mathcal{V}, \mathcal{C}, \|\cdot\|)$ may be induced by multiple gauges; i.e. there may exist different matrix gauges $\nu$ and $\omega$ on $\mathcal{V}$ which satisfy $\mathcal{C} = \mathcal{C}^{\nu} = \mathcal{C}^{\omega}$ and $\| \cdot \| = \|\cdot\|^{\nu} = \|\cdot\|^{\omega}$. The following definition will show us how to define a canonical matrix gauge on a given normal matrix ordered $*$-operator space which induces its matrix ordering and matrix norm.

\begin{definition}
Suppose that $(\mathcal{V}, \mathcal{C}, \|\cdot\|)$ is a normal matrix ordered $*$-operator space. For each $x \in M_n(\mathcal{V})_h$, we define
\[ \nu_n^{max}(x) := \inf \{ \| x + p \|_n : p \in \mathcal{C}_n \}. \]
We let $\nu^{max}$ denote the sequence $\{\nu_n^{max}\}_{n=1}^\infty$.
\end{definition}

\begin{example}
Let $\mathcal{V}$ be the one-dimensional self-adjoint operator space from Example \ref{example: non-gauge isometric spaces}, and let $A$ be the corresponding matrix. Then
$\nu_1^{max}(A) = \nu_1^{max}(-A)=\|\pm A\|=1$, since the positive cone of $\mathcal{V}$ consists of only the zero vector.
\end{example}

\begin{theorem} \label{thm: max gauge}
Suppose that $(\mathcal{V}, \mathcal{C}, \| \cdot \|)$ is a normal matrix ordered $*$-vector space. Then $\nu^{max}$ is a proper matrix gauge on $\mathcal{V}$, $\mathcal{C} = \mathcal{C}^{\nu^{max}}$ and $\| \cdot \| = \|\cdot\|_{\nu^{\max}}$. If $\nu$ is another proper matrix gauge satisfying $\mathcal{C} = \mathcal{C}^{\nu}$ and $\|\cdot \| = \|\cdot\|^{\nu}$, then $\nu_n(x) \leq \nu_n^{max}(x)$ for all $x \in M_n(\mathcal{V})$.
\end{theorem}

\begin{proof}
First, let $x,y \in M_n(\mathcal{V})_h$, and suppose that $p, q \in \mathcal{C}_n$. Then
\[ \|x + y + p + q \|_n \leq \|x + p\|_n + \|y + q\|_n. \]
It follows that $\nu_n^{max}(x+y) \leq \nu_n^{max}(x) + \nu_n^{max}(y)$. Next, assume $\alpha \in M_{n,k}$. Then for any $p \in \mathbb{C}_n$,
\[ \|\alpha^*x \alpha + \alpha^*p\alpha\|_k \leq \|\alpha\|^2 \|x + p \|_n. \]
Hence $\nu_k^{max}(\alpha^* x \alpha) \leq \|\alpha\|^2 \nu_n^{max}(x)$. Finally, let $z \in M_m(\mathcal{V})_h$. Then for any $p \in \mathbb{C}_n$ and $q \in \mathbb{C}_m$, we have
\[ \|x \oplus y + p \oplus q\|_{n+m} = \max(\|x + p\|_n, \|y + q\|_m). \]
Thus, $\nu_{n+m}^{max}(x \oplus y) \leq \max( \nu_n^{max}(x), \nu_m^{max}(y))$. For the other inequality, Suppose that $R \in \mathcal{C}_{n+m}.$ Let 
\[ p = \begin{pmatrix} I_n & 0_{n,m} \end{pmatrix} R \begin{pmatrix} I_n \\ 0_{m,n} \end{pmatrix} \quad \text{ and } q = \begin{pmatrix} 0_{m,n} & I_m \end{pmatrix} R \begin{pmatrix} 0_{n,m} \\ I_m \end{pmatrix}. \]
Then
\[ \|x + p\|_n = \| \begin{pmatrix} I_n & 0_{n,m} \end{pmatrix} (x \oplus y + R) \begin{pmatrix} I_n \\ 0_{m,n} \end{pmatrix} \|_n \leq \|x \oplus y + R \|_{n+m} \]
and
\[ \|y + q\|_m = \| \begin{pmatrix} 0_{m,n} & I_m \end{pmatrix} (x \oplus y + R) \begin{pmatrix} 0_{n,m} \\ I_m \end{pmatrix} \|_m \leq \|x \oplus y + R \|_{n+m}. \]
Thus, $\max(\nu_n^{max}(x), \nu_m^{max}(y)) \leq \nu_{n+m}^{max}(x \oplus y)$. Therefore $\nu^{max}$ is a matrix gauge. We will show that $\nu^{max}$ is proper below.

Next, we check that $\mathcal{C}^{\nu^{max}} = \mathcal{C}$. Observe that
\[ \nu_n^{max}(-x) = \inf\{ \|-x+p\|_n : p \in \mathcal{C}_n\} = \inf \{ \|x - p\|_n : p \in \mathcal{C}_n\}. \]
Thus, $\nu_n(-x)$ is precisely the distance between $x$ and the closed cone $\mathcal{C}_n$. Since $\mathcal{C}_n$ is closed, we see that $\nu_n(-x)=0$ if and only if $x \in \mathcal{C}_n$, proving the claim. It follows that $\nu^{max}$ is proper, since $\nu_n^{max}(x) = \nu_n^{max}(-x) = 0$ implies that $\pm x \in \mathcal{C}_n$, and $\pm x \in \mathcal{C}_n$ implies that $x = 0$, since $\mathcal{C}_n$ is a proper cone.

Next, we check that $\|\cdot\| = \|\cdot\|^{\nu^{max}}$. First, since $x \mapsto x^*$ is completely isometric with respect to $\|\cdot\|$, we have 
\[ \|x\|_n = \| \begin{pmatrix} 0 & x \\ x^* & 0 \end{pmatrix} \|_{2n} \]
for every $x \in M_n(\mathcal{V})$. Let $x \in M_n(\mathcal{V})$. Then 
\[ \nu_{2n} \begin{pmatrix} 0 & \pm x \\ \pm x^* & 0 \end{pmatrix} \leq \| \begin{pmatrix} 0 & \pm x \\ \pm x^* & 0 \end{pmatrix} + 0 \|_{2n} = \|x\|_n \]
since $0 \in \mathcal{C}_{2n}$. Hence, $\|x\|_{\nu^{max},n} \leq \|x\|_n$. Furthermore, since $\mathcal{V}$ is normal, and since
\[ \begin{pmatrix} 0 & x \\  x^* & 0 \end{pmatrix} - Q \leq \begin{pmatrix} 0 & x \\  x^* & 0 \end{pmatrix} \leq \begin{pmatrix} 0 &  x \\  x^* & 0 \end{pmatrix} + P \]
for all $Q,P \in \mathcal{C}_{2n}$, we have
\[ \|x\|_n \leq \max( \| \begin{pmatrix} 0 & - x \\ - x^* & 0 \end{pmatrix} + Q \|_{2n}, \| \begin{pmatrix} 0 & x \\ x^* & 0 \end{pmatrix} + P \|_{2n}). \]
It follows that $\|x\|^{\nu^{max}} _n \geq \|x\|_n$.

Finally, suppose that $\nu$ is another matrix gauge satisfying $\mathcal{C}^{\nu} = \mathcal{C}$ and $\|\cdot\|^{\nu} = \|\cdot\|$. Let $x \in M_n(\mathcal{V})_h$ and let $p \in \mathcal{C}_n$. Then
\begin{eqnarray}
\nu_n(x) & = & \nu_n(x + p - p) \nonumber \\
& \leq & \nu_n(x + p) + \nu_n(-p) \nonumber \\
& = & \nu_n(x + p) \nonumber \\
& \leq & \|x + p\|_n. \nonumber
\end{eqnarray}
It follows that $\nu_n(x) \leq \nu_n^{max}(x)$.
\end{proof}

\begin{example} \label{ex: order guage is maximal}
Let $\mathcal{S}$ be an operator system with unit $e$. Then $\nu^{max} = \nu^{e}$. To see this, let $x \in M_n(\mathcal{S})_h$, and suppose that $a > \nu_n^e(x)$ and $b > \nu_n^e(-x)$. Then $ae - x \geq 0$ and $x + be \geq 0$. Now
\[ x = \frac{a}{a+b}(x+be) - \frac{b}{a+b}(ae - x). \]
Moreover,
\[ \nu_n^e(\frac{a}{a+b}(x+be)) \leq \frac{a}{a+b} (\nu_n^e(x) + b) \leq \frac{a}{a+b}(a+b) = a \]
and
\[ \nu_n^e(\frac{b}{a+b}(ae-x)) \leq \frac{b}{a+b} (a + \nu_n^e(-x)) \leq \frac{b}{a+b}(a+b) = b. \]
It follows that 
\[ \nu_n^{max}(x) \leq \|x + \frac{b}{a+b}(ae-x)\|_n^e = \|\frac{a}{a+b}(x+be))\|_n^e = \frac{a}{a+b}\nu_n^e(x+be) \leq a. \]
Since this holds for every $a > \nu_n^e(x)$, we conclude that $\nu_n^e(x) \geq \nu_n^{max}(x)$. The other inequality follows from Theorem \ref{thm: max gauge}.
\end{example}

As a corollary of Theorem \ref{thm: max gauge} and Corollary \ref{cor: Representations of matrix gauge spaces}, we obtain a proof of Werner's claim that normal matrix ordered operator spaces admit completely isometric complete order embeddings into $B(H)$.

\begin{corollary} \label{cor: normal OVS}
Let $(\mathcal{V}, \mathcal{C}, \|\cdot\|)$ be a normal matrix ordered $*$-operator space. Then there exists a Hilbert space $H$ and a completely isometric complete order embedding $\phi: \mathcal{V} \to B(H)$.
\end{corollary}

\begin{proof}
By Theorem \ref{thm: max gauge}, we have $\mathcal{C}=\mathcal{C}^{\nu^{max}}$ and $\|\cdot\| = \|\cdot\|^{\nu^{max}}$. By Corollary \ref{cor: Representations of matrix gauge spaces}, there exists a Hilbert space $H$ and a completely gauge isometric embedding $\phi: \mathcal{V} \to B(H)$. By Proposition \ref{prop: gauge contractive is cp and cc}, $\phi$ is a completely isometric complete order embedding.
\end{proof}


In Example \ref{example: non-gauge isometric spaces}, we showed that maps which are completely contractive and completely positive on the operator space $\mathcal{V}$ do not necessarily have extensions that are both completely contractive and completely positive. In particular, a completely positive completely contractive linear functional on a non-unital operator space may not extend to a state on the unitization of the given space. We also showed that the restriction of $\nu^e$ on $M_2$ to $\mathcal{V}$ differs from the maximal gauge on $\mathcal{V}$. In the following results, we will show that these observations are not coincidental.

\begin{theorem} \label{thm: cpcc equals max cgc}
Let $(\mathcal{V}, \mathcal{C}, \|\cdot\|)$ be a normal matrix ordered $*$-operator space and let $(\mathcal{W}, \omega)$ be a matrix gauge $*$-vector space. Then a self-adjoint linear map $\phi: \mathcal{V} \to \mathcal{W}$ is completely positive and completely contractive if and only if it is completely gauge contractive as a map from $(\mathcal{V}, \nu^{max})$ to $(\mathcal{W}, \omega)$.
\end{theorem}

\begin{proof}
First suppose that $\phi$ is completely positive and completely contractive. Suppose that $x \in M_n(\mathcal{V})$ and $x=x^*$. Then for every positive $p \in M_n(\mathcal{V})$,
\begin{eqnarray}
\omega_n(\phi^{(n)}(x)) & = & \omega_n(\phi^{(n)}(x) + \phi^{(n)}(p) - \phi^{(n)}(p)) \nonumber \\
& \leq & \omega_n(\phi^{(n)}(x) + \phi^{(n)}(p)) + \omega_n(-\phi^{(n)}(p)) \nonumber \\
& = & \omega_n(\phi^{(n)}(x) + \phi^{(n)}(p)) \nonumber \\
& \leq & \|\phi_n(x + p)\|_n \nonumber \\
& \leq & \|x + p\|_n.
\end{eqnarray}
Taking an infimum over all positive $p \in M_n(\mathcal{V})$ gives $\omega_n(\phi^{(n)}(x)) \leq \nu_n^{max}(x)$. So $\phi$ is completely gauge contractive.

Now assume that $\phi$ is completely gauge contractive. Then $\phi$ is completely positive and completely contractive by Proposition \ref{prop: gauge contractive is cp and cc}.
\end{proof}

We are now prepared to address the question of when completely contractive completely positive maps have completely contractive completely positive extensions. 

\begin{theorem} \label{thm: extension of cpcc}
Let $(\mathcal{V}, \mathcal{C}, \|\cdot\|)$ be a normal matrix ordered operator space, and let $\mathcal{W}$ be a self-adjoint subspace. Let $\phi: \mathcal{W} \to B(H)$ be a completely positive completely contractive map. Then there exists a completely positive completely contractive extension $\tilde{\phi}: \mathcal{V} \to B(H)$ if and only if $\phi$ is completely gauge contractive with respect to the restriction of the maximal gauge $\nu^{max}$ on $\mathcal{V}$ to $\mathcal{W}$.
\end{theorem}

\begin{proof}
Let $\omega$ denote the restriction of $\nu^{max}$ on $\mathcal{V}$ to $\mathcal{W}$. Then $\mathcal{W}$ induces the norm and order structure on $\mathcal{W}$, since it induces the norm and order structure on $\mathcal{V}$ and the norm and order structure on $\mathcal{W}$ is inherited from $\mathcal{V}$.

Suppose that $\phi$ is completely gauge contractive with respect to $\omega$. Then $\phi$ is completely positive and completely contractive, by Proposition \ref{prop: gauge contractive is cp and cc}. By Theorem \ref{thm: gauge Arveson}, $\phi$ has a completely gauge contractive extension $\widetilde{\phi}: \mathcal{V} \to B(H)$. Since $\nu^{max}$ induces the norm and order structure on $\mathcal{V}$, $\widetilde{\phi}$ is completely positive and completely contractive, by Proposition \ref{prop: gauge contractive is cp and cc}.

On the other hand, suppose that $\phi$ admits an extension $\tilde{\phi}: \mathcal{V} \to B(H)$ which is completely positive and completely contractive. Then $\tilde{\phi}$ is competely gauge contractive with respect to $\nu^{max}$ on $\mathcal{V}$, by Theorem \ref{thm: cpcc equals max cgc}. It follows that the restriction $\phi$ to $\mathcal{W}$ is completely gauge contractive with respect to $\omega$.
\end{proof}

The next definition will allow us to characterize when every completely positive completely contractive map admits an extension.

\begin{definition} \label{defn: gauge maximal}
Let $(\mathcal{V}, \mathcal{C}, \|\cdot\|)$ be a normal matrix ordered operator space, and let $\mathcal{W}$ be a self-adjoint subspace. We say that $\mathcal{W}$ is \textbf{gauge maximal} if the restriction of $\nu^{max}$ on $\mathcal{V}$ to $\mathcal{W}$ coincides with the maximal gauge on $\mathcal{W}$. In other words, we ask that for each $x \in M_n(\mathcal{W})_h$,
\[ \inf \{ \|x + p\| : p \in M_n(\mathcal{W}) \cap \mathcal{C}_n \} = \nu^{max}(x) = \inf \{ \|x + p\| : p \in \mathcal{C}_n \}. \]
\end{definition}

We will see an example of a subspace which is not gauge maximal after the next result.  However, because the infimum on the left is taken over a smaller set than the infimum on the right in Definition \ref{defn: gauge maximal}, the restriction of $\nu^{max}$ to a subspace $\mathcal{W}$ should be expected to be smaller than the maximal gauge on $\mathcal{W}$ in general. The following corollary is immediate from Theorem \ref{thm: extension of cpcc}

\begin{corollary} \label{cor: extensions of cpcc and gague maimal}
Let $(\mathcal{V}, \mathcal{C}, \|\cdot\|)$ be a normal matrix ordered operator space, and let $\mathcal{W}$ be a self-adjoint subspace. Then every completely positive completely contractive map $\varphi: \mathcal{W} \to B(H)$ admits a completely positive completely contractive extension $\widetilde{\phi}: \mathcal{V} \to B(H)$ if and only if $\mathcal{W}$ is a gauge maximal subspace of $\mathcal{V}$.
\end{corollary}

\begin{remark}
\emph{Consider the category whose objects are normal matrix ordered operator spaces and whose morphisms are completely positive completely contractive maps. Theorem \ref{thm: extension of cpcc} implies that $B(H)$ is not injective in this category. In fact, this holds regardless of $\dim(H)$ and even in the case when $\dim(H)=1$, as illustrated in Example \ref{example: non-gauge isometric spaces}. In fact, if $\varphi: \mathcal{V} \to \mathbb{C}$ is the linear functional described in Example \ref{example: non-gauge isometric spaces}, any contractive extension of $\varphi$ to $\mathcal{V}'$ mapping the unit to a positive value will fail to be positive. This failure is a consequence of the fact that $\mathcal{V}$ is not gauge maximal in $\mathcal{V}'$.}
\end{remark}

In spite of Theorem \ref{thm: extension of cpcc}, there are many cases where completely positive completely contractive maps are known to always have completely positive completely contractive extensions, such as unital inclusions of operator systems or inclusions of (possibly non-unital) C*-algebras. Each of these situations involves the inclusion of a gauge maximal subspace. We can check these cases directly.

\begin{theorem} \label{thm: opSys gauge maximal}
Let $\mathcal{T} \subseteq \mathcal{S}$ be a unital inclusion of operator systems. Then $\mathcal{T}$ is gauge maximal.
\end{theorem}

\begin{proof}
This follows from Example \ref{ex: order guage is maximal}.
\end{proof}

\begin{theorem} \label{thm: Cstar gauge maximal}
Let $\mathcal{B}$ be a (possibly non-unital) C*-algebra, and let $\mathcal{A}$ be a C*-subalgebra. Then $\mathcal{A}$ is a gauge maximal subspace. Hence every completely positive completely contractive map $\phi: \mathcal{A} \to B(H)$ admits a completely positive completely contractive extension $\widetilde{\phi}: \mathcal{B} \to B(H)$.
\end{theorem}

\begin{proof}
If $\mathcal{A}$ is unital, then so is $\mathcal{B}$ and the result follows from Theorem \ref{thm: opSys gauge maximal}. So assume this is not the case. Let $x \in M_n(\mathcal{A})$ be self-adjoint. Let $C^*(x)$ denote the C*-subalgebra of $M_n(\mathcal{A})$ generated by $x$. Since $M_n(\mathcal{A})$ is non-unital, $C^*(x)$ is isometrically $*$-isomorphic to $C_0(\sigma(x) \setminus \{0\})$, where $\sigma(x)$ denotes the spectrum of $x$ in $\mathbb{R}$ and $C_0(\sigma(x) \setminus \{0\})$ denotes the set of continuous functions on $\sigma(x)$ vanishing at $0$ (see Corollary I.2.6 of \cite{DavidsonCstarBook}), via the Gelfand transform $\Gamma$ (mapping $x$ to the function $f(y)=y$). Define
\[ f_+(y) = \begin{cases} y & y \geq 0 \\ 0 & \text{else} \end{cases} \qquad \text{and} \qquad f_-(y) = \begin{cases} 0 & y \geq 0 \\ y & \text{else} \end{cases} \]
and observer that $f_+, f_- \in C_0(\sigma(x) \setminus \{0\})$. Morever, $\Gamma(x_+) = f_+$, $\Gamma(x_-) = f_-$, and
\[ \nu_n^{max}(x) \leq \|x + x_-\| = \|f + f_-\| = \|f_+\| = \|x_+\|. \]
Since $x_- \in M_n(\mathcal{A})$, it follows that $\mathcal{A}$ is a gauge maximal subspace of $\mathcal{B}$.
\end{proof}

We remark that Theorem \ref{thm: Cstar gauge maximal}, together with Corollary \ref{cor: Unique unital extension op sys}, implies Theorem 2.2.1 of \cite{Ozawa2008CAlgebrasAF}, whose proof relied on properties of the Banach space double dual of a C*-algebra.

\section{Kernels of completely positive maps} \label{sec: kernels and quotients}

We continue our study of ordered operator spaces with some results
regarding quotients. To motivate this topic, we first consider the question of defining a quotient on an ordered vector space.

Let $(V,V_+)$ be an ordered vector space, i.e. $V$ is a real vector space and $V_+ \subseteq V$ is a proper cone. We call a subspace $J \subseteq V$ an \textbf{order ideal} provided that whenever $0 \leq x \leq y$ for some $x \in V$ and $y \in J$, it must be that $x \in J$. It is easy to see that if $\phi$ is a positive linear map from $V$ to some other ordered vector space $W$, then $\ker(\phi)$ is an order ideal. In fact, every order ideal can be realized as the kernel of a positive map. This is done by regarding the algebraic quotient $V / J$ of $V$ by an order ideal $J$ as an ordered vector space by setting $(V/J)_+ = \{p + J: p \in V_+\}$. It is straightforward to check that $(V/J)_+$ is a proper positive cone. Hence, the quotient map $\pi: V \rightarrow (V/J)_+$ is positive, and $J = \ker(\pi)$.

When $V$ is a normed ordered space, the situation becomes less straightforward. It is easy to check that the kernel of any bounded (i.e., continuous) positive map is a norm closed order ideal in $V$. However, not every norm closed order ideal in $V$ is the kernel of a bounded positive map. For example, Paulsen and Tomforde show in \cite{PaulsenTomforde} that the span of the matrix unit $E_{1,1}$ in $(M_2)_{h}$ is not the kernel of a bounded positive linear map, despite being a closed order ideal. This problem is closely related to the problem of defining quotients of ordered vector spaces. If $J$ is a non-unital closed order ideal in an AOU space $V$, then the ordered space $(V/J)$ is an order unit space with order unit $e + J$, but this order unit may not be Archimedean. Paulsen and Tomforde developed an Archimedeanization process which converts an order unit space to an AOU space by identifying a subspace of infinitesimals and showing that the quotient by this space is an AOU space. Hence, one can define quotients of AOU spaces by non-unital order ideals in terms of the Archimedeanization process. These technicalities persist in the theory of operator systems as well. For example, Section 3 of \cite{Kavruk} is largely devoted to identifying which order ideals in an operator system can be realized as kernels of completely positive maps, and the Archimedeanization process plays an important role in their analysis.

In the following, we will characterize those subspaces of matrix gauge $*$-vector spaces which arise as kernels of completely gauge contractive maps. These kernels arise as \textit{gauge ideals}, which we define next. These results also characterize kernels of completely positive maps on operator systems, since an operator system has a matrix gauge structure, and so we will point out several applications to operator systems below.

\begin{definition}
Let $(\mathcal{V}, \nu)$ be a matrix gauge $*$-vector space.
\begin{enumerate}
    \item We call a self-adjoint subspace $\mathcal{J} \subseteq \mathcal{V}$ a \textbf{gauge ideal} provided that whenever there exist sequences $\{a_k\}_{k=1}^\infty, \{b_k\}_{k=1}^\infty \subseteq \mathcal{J}_h$ and an element $x \in \mathcal{V}_h$ satisfying 
    \[ \nu_1(x - a_k) \to 0 \quad \text{and} \quad \nu_1(b_k - x) \to 0 \]
    then $x \in \mathcal{J}$.
    \item We call a self-adjoint subspace $\mathcal{J} \subseteq \mathcal{V}$ a \textbf{complete gauge ideal} provided that, for every $n \in \mathbb{N}$, whenever there exist sequences $\{a_k\}_{k=1}^\infty, \{b_k\}_{k=1}^\infty \subseteq M_n(\mathcal{J})_h$ and an element $x \in M_n(\mathcal{V})_h$ satisfying 
    \[ \nu_n(x - a_k) \to 0 \quad \text{and} \quad \nu_n(b_k - x) \to 0 \]
    then $x \in M_n(\mathcal{J})$.
\end{enumerate}
\end{definition}

It is clear that every complete gauge ideal is a gauge ideal. It turns out that a gauge ideal is automatically a complete gauge ideal. This is proven in the following lemma.

\begin{lemma} \label{elinfQuot} Let $(\mathcal{V},\nu)$ be an matrix gauge $*$-vector space. Suppose that $\mathcal{J} \subseteq \mathcal{V}$ is a gauge ideal. Then $\mathcal{J}$ is a complete gauge ideal.
\end{lemma}

\begin{proof} Let $C = (c_{i,j}) \in M_n(\mathcal{V})_h$, and suppose that there exist sequences $\{A_n\},\{B_n\} \subset M_n(\mathcal{J})_h$ such that $\nu_n(A_n + C) \to 0$ and $\nu_n(B_n - C) \rightarrow 0$. Pick an integer $j \in \{1,\dots, n\}$, and let $e_j$ be the $1 \times n$ matrix with a 1 in the $j$th entry and zeroes elsewhere. Then $\nu_1(e_j(A_n + C)e_j^*) \leq \nu_n(A_n + C) \rightarrow 0$ and $\nu_1(e_j(B_n - C)e_j^*) \leq \nu_n(B_n - C) \rightarrow 0$. Since $\mathcal{J}$ is a gauge ideal of $\mathcal{V}$, we see that $c_{j,j} \in \mathcal{J}$. Hence, all the diagonal entries of $C$ are in $\mathcal{J}$. Also, for each pair of integers $k,j \in \{1,\dots,n\}$, the sequences $\nu_1((e_k + e_j)(A_n + C)(e_k + e_j)^*)$ and $\nu_1((e_k + e_j)(B_n - C)(e_k + e_j)^*)$ converge to zero. It follows that that $c_{k,j} + c_{k,j}^* \in \mathcal{J}$. Likewise, $\nu_1((e_k + ie_j)(A_n + C)(e_k + ie_j)^*)$ and $\nu_1((e_k + ie_j)(B_n - C)(e_k + ie_j)^*)$ converge to zero, and hence, $i(c_{k,j} - c_{k,j}^*) \in \mathcal{J}$. It follows that $c_{k,j} \in \mathcal{J}$ for each $k,j \in \{1, \dots, n\}$, and hence, $C \in M_n(\mathcal{J})_h$. So $\mathcal{J}$ is a complete gauge ideal of $\mathcal{V}$.
\end{proof}

Let $\mathcal{J}$ be a gauge ideal in $\mathcal{V}$, and suppose that $0 \leq x \leq y$ for some $x \in \mathcal{V}$ and $y \in \mathcal{J}$. Then $\nu_1(0-x)=\nu_1(x-y) = 0$. By regarding $\{0\}$ and $\{y\}$ as constant sequences, we see that $x \in \mathcal{J}$. So $\mathcal{J}$ is automatically an order ideal. We will show subsequently that $\mathcal{V} / \mathcal{J}$ may be regarded as a matrix gauge $*$-vector space. To do so, we must define a matrix gauge on the vector space $\mathcal{V} / \mathcal{J}$.

\begin{definition}
Let $(\mathcal{V},\nu)$ be a matrix gauge $*$-vector space, and suppose that $\mathcal{J} \subseteq \mathcal{V}$ is a gauge ideal. For each $n \in \mathbb{N}$, we define a function $q_n: M_n(\mathcal{V})_h / M_n(\mathcal{J})_h \to [0,\infty)$ by
\[ q_n(A + M_n(\mathcal{J})_h) = \inf \{ \nu_n(B) : B \in A + M_n(\mathcal{J})_h \}. \]
We let $q^{\mathcal{J}}$ (or simply $q$ when $\mathcal{J}$ is clear from context) denote the sequence $\{q_n\}_{n=1}^\infty$.
\end{definition}

\begin{theorem} Let $(\mathcal{V},\nu)$ be an matrix gauge $*$-vector space, and let
$\mathcal{J}$ be a gauge ideal of $\mathcal{V}$. Then $q^{\mathcal{J}}$ is a matrix gauge on $\mathcal{V} / \mathcal{J}$. \end{theorem}

\begin{proof} We first check that $q := q^{\mathcal{J}}$ is proper. Suppose that $q_n(\pm x + M_n(\mathcal{J})_h) = 0$ for some $x \in M_n(\mathcal{V})_h$. Then there exist sequences $\{a_n\}, \{b_n\} \subseteq M_n(\mathcal{J})_h$ such that 
\[ \nu_n(x + a_n) \to 0 \quad \text{and} \quad \nu_n(b_n - x) \to 0. \]
Therefore $x \in M_n(\mathcal{J})_h$ and hence $x + M_n(\mathcal{J})_h = M_n(\mathcal{J})_h$. So $q$ is proper.

Let $A \in M_n(\mathcal{V})_h$ and $X \in M_{n,k}$. 
Then
\begin{eqnarray} q_k(X^*AX + M_k(\mathcal{J})_h) & = & \inf\{ \nu_k(B) : B
\in X^*AX + M_k(\mathcal{J})_h \} \nonumber \\ & \leq & \inf \{
\nu_k(X^*CX) : C \in A + M_n(\mathcal{J})_h \} \nonumber \\ & \leq &
\|X\|^2 q_n(A + M_n(\mathcal{J})_h). \nonumber \end{eqnarray} Now, let $A
\in M_n(\mathcal{V})_h$ and $B \in M_k(\mathcal{V})_h$. For any $R \in M_{n,k}(\mathcal{V})$ we have
\[ \nu_n(A) = \nu_n( P_1 
\begin{pmatrix}
A & R \\
R^*  & B \end{pmatrix} P_1^* ) \leq \|P_1\|^2  \nu_{n+k} 
\begin{pmatrix}
A & R \\
R^*  & B \end{pmatrix} = \nu_{n+k} 
\begin{pmatrix}
A & R \\
R^*  & B \end{pmatrix}\]
and
\[ \nu_k(B) = \nu_k( P_2 
\begin{pmatrix}
A & R \\
R^*  & B \end{pmatrix} P_2^* ) \leq \|P_2\|^2  \nu_{n+k} 
\begin{pmatrix}
A & R \\
R^*  & B \end{pmatrix} = \nu_{n+k} 
\begin{pmatrix}
A & R \\
R^*  & B \end{pmatrix}, \]
where $P_1 = (I_n \quad 0_{n,k}) \in M_{n,n+k}, P_2 = (0_{k,n} \quad I_k) \in M_{k,n+k}$. Thus \[ \nu_{n+k} 
\begin{pmatrix}
A & R \\
R^*  & B \end{pmatrix}  \geq \nu_{n+k}
\begin{pmatrix}
A & 0 \\
0 & B \end{pmatrix}  \] for any $R \in M_{n,k}(\mathcal{J})$, since the
right hand side is equal to $\max\{\nu_n(A),\nu_k(B)\}$. Hence,
\begin{eqnarray} q_{n+k}(A \oplus B + M_{n+k}(\mathcal{J})_h) & = & \inf
\{\nu_{n+k}(A' \oplus B') : A' \in A + M_n(\mathcal{J})_h, B' \in B +
M_k(\mathcal{J})_h \} \nonumber \\ & = & \max\{q_n(A + M_{n}(\mathcal{J})_h),q_k(B
+ M_{k}(\mathcal{J})_h)\}. \nonumber
\end{eqnarray}
Therefore $q$ is a proper matrix gauge on $\mathcal{V}/\mathcal{J}$. \end{proof}

\begin{corollary} \label{cor: cgc kernel is gauge ideal}
Let $(\mathcal{V}, \nu)$ be a matrix gauge $*$-vector space. Then a subspace $\mathcal{J} \subseteq \mathcal{V}$ is the kernel of a completely gauge contractive map if and only if it is a gauge ideal.
\end{corollary}

\begin{proof}
First suppose that $\mathcal{J}$ is a gauge ideal. Then $\mathcal{J}$ is the kernel of the quotient map $\pi: \mathcal{V} \to \mathcal{V} / \mathcal{J}$. The map $\pi$ is completely gauge contractive since for every $x \in M_n(\mathcal{V})_h$ we have $q_n(\pi(x)) = q_n(x + M_n(\mathcal{J})_h) \leq \nu_n(x+0)=\nu_n(x)$.

Conversely, suppose that $\mathcal{J} = \ker(\phi)$ where $\phi: \mathcal{V} \to \mathcal{W}$ is completely gauge contractive, where $(\mathcal{W}, \omega)$ is another matrix gauge $*$-vector space. Suppose that there exist sequences $\{a_n\}, \{b_n\} \in \mathcal{J}$ such that $\nu_1(a_n - x) \to 0$ and $\nu_1(x - b_n) \to 0$ for some $x \in \mathcal{V}$. Then
\[ \omega_1( -\phi(x) ) = \omega_1( \phi(a_n - x) ) \leq \nu_1(a_n - x) \to 0 \]
and
\[ \omega_1( \phi(x) ) = \omega_1( \phi(x - b_n) ) \leq \nu_1(x - b_n) \to 0. \]
Therefore $\omega_1(\pm \phi(x)) = 0$ and hence $\phi(x) = 0$, since $\omega$ is proper. So $x \in \ker(\phi) = \mathcal{J}$. We conclude that $\mathcal{J}$ is a gauge ideal.
\end{proof}

Since every operator system is a matrix gauge $*$-vector space with respect to the gauge $\nu^e$, we have also characterized the kernels of completely positive maps.

\begin{corollary} \label{cor: kernels of cp maps}
Let $\mathcal{J}$ be a self-adjoint subspace of an operator system $\mathcal{S}$. Then $\mathcal{J}$ is the kernel of a completely positive map if and only if it is a gauge ideal of the matrix gauge $*$-vector space $(\mathcal{S}, \nu^e)$. Moreover, the quotient $\mathcal{S}/\mathcal{J}$ is an operator system with unit $\hat{e} := e + \mathcal{J}$, and $q^{\mathcal{J}} = \nu^{\hat{e}}$.
\end{corollary}

\begin{proof}
For the first statement, if $\phi$ is completely positive, then it is a scalar multiple of a completely positive completely contractive map. Hence $\mathcal{J}$ is the kernel of a completely positive map if and only if it is the the kernel of a completely positive completely contractive map. By Theorem \ref{thm: cpcc equals max cgc}, the completely positive completely contractive maps coincide with the completely gauge contractive maps with respect to the maximal gauge $\nu^{max}$. However, $\nu^e = \nu^{max}$ by Example \ref{ex: order guage is maximal}. The result follows from Corollary \ref{cor: cgc kernel is gauge ideal}.

For the second statement, define for each $n \in \mathbb{N}$ the set $C_n(\mathcal{S}/\mathcal{J})$ to be the cosets $x + M_n(\mathcal{J})_h$ such that for every $t > 0$ there exists $y_t \in M_n(\mathcal{J})_h$ such that $0 \leq x + y_t + t I_n \otimes e$. By Proposition 3.4 of \cite{Kavruk}, $(\mathcal{S} / \mathcal{J}, \{C_n(\mathcal{S} / \mathcal{J})\}, e + \mathcal{J})$ is an operator system. Suppose that $q_n(-x + M_n(\mathcal{J})_h) = 0$. Then for every $t > 0$, there exists $y_t \in M_n(\mathcal{J})_h$ such that $\nu_n^e(-x + y_t) < t/2$. Hence, for every $t > 0$, there exists $y_t \in M_n(\mathcal{J})_h$ such that $0 \leq x - y_t + t I_n \otimes e$. It follows that $x \in C_n(\mathcal{S} / \mathcal{J})$. On the other hand, suppose that $x \in C_n(\mathcal{S} / \mathcal{J})$. Then for every $t > 0$, there exists $y_t \in M_n(\mathcal{J})_h$ such that $0 \leq x + y_t + t I_n \otimes e$ and hence $\nu_n^e(-x - y_t) < t$. Since this holds for all $t > 0$, $q_n(-x + M_n(\mathcal{J})_h) = 0$. It follows that $\mathcal{C}_n^{q} = C_n(\mathcal{S} / \mathcal{J})$ for every $n \in \mathbb{N}$. Thus $(\mathcal{S}/\mathcal{J}, \mathcal{C}^q, \hat{e})$ is an operator system. It remains to verify that $\nu^{\hat{e}} = q$. This follows from checking that, for each $x \in M_n(\mathcal{S})_h$, the following statements are equivalent:
\begin{itemize}
    \item $\nu_n^{\hat{e}}(x + M_n(\mathcal{J})_h) \leq t$.
    \item $x + M_n(\mathcal{J})_h \leq te + M_n(\mathcal{J})_h$.
    \item $q(x - te + M_n(\mathcal{J})_h) = 0$.
    \item For every $\epsilon > 0$, there exists $y_{\epsilon} \in M_n(\mathcal{J})_h$ such that $x - te + y_{\epsilon} \leq \epsilon I_n \otimes e$.
    \item For every $\epsilon > 0$ there exists $y_{\epsilon} \in M_n(\mathcal{J})_h$ such that $x + y_{\epsilon} \leq (t+\epsilon) I_n \otimes e$.
    \item $q(x+M_n(\mathcal{J})_h) \leq t$.
\end{itemize}
Since $\nu_n^{\hat{e}}(x + M_n(\mathcal{J})_h) \leq t$ if and only if $q(x+M_n(\mathcal{J})_h) \leq t$, $\nu^{\hat{e}} = q$.
\end{proof}

As a demonstration of Corollary \ref{cor: kernels of cp maps}, we show directly that the linear span of $E_{1,1}$ is not the kernel of a completely positive map on $M_2$.

\begin{example} \emph{Consider $M_2$ as a matrix gauge $*$-vector space equipped with the gauge $\nu^e$. Recall that a matrix $A \in M_2$ is positive if and only if the diagonal entries of $A$ are positive and $\det(A) \geq 0$. Fix $\lambda \in \mathbb{C}$ with $|\lambda|=1$. Then for all $n \in \mathbb{N}$,}
\[ \det 
\begin{pmatrix}
n + \frac{1}{n} & \pm \lambda \\
\pm \overline{\lambda}  & \frac{1}{n} \end{pmatrix} = \frac{1}{n^2} > 0. \]
\emph{Hence},
\[
\begin{pmatrix}
-n  & 0 \\
0  & 0 \end{pmatrix} \pm 
\begin{pmatrix}
0 & \lambda \\
\overline{\lambda}  & 0 \end{pmatrix} \leq 
\begin{pmatrix}
\frac{1}{n} & 0 \\
0 & \frac{1}{n} \end{pmatrix} \] \emph{for each $n \in
\mathbb{N}$. So}
\[ \nu_1^e \left(
\begin{pmatrix}
-n  & 0 \\
0  & 0 \end{pmatrix} -
\begin{pmatrix}
0 & \lambda \\
\overline{\lambda}  & 0 \end{pmatrix} \right) \rightarrow 0. \]
\emph{Similarly,} \[ \nu_1^e \left( 
\begin{pmatrix}
-n  & 0 \\
0  & 0 \end{pmatrix}  + 
\begin{pmatrix}
0 & \lambda \\
\overline{\lambda}  & 0 \end{pmatrix}  \right) \rightarrow 0. \]
\emph{It follows that the span of $E_{1,1}$ is not a gauge ideal, since}
\[ 
\begin{pmatrix}
0 & \lambda \\
\overline{\lambda}  & 0 \end{pmatrix}  \]
\emph{is not an element of the span of $E_{1,1}$.}
\end{example}

We conclude this section with one last application of Corollary \ref{cor: kernels of cp maps}. Using gauge quotients, we can characterize the matrix gauge $*$-vector spaces $(\mathcal{V}, \nu)$ which are completely gauge isometric to an operator system with its canonical gauge $\nu^e$. To do so, we must first introduce complete gauge quotient maps.

\begin{definition}
Let $(\mathcal{V}, \nu)$ and $(\mathcal{W}, \omega)$ be matrix gauge $*$-vector spaces, and let $\pi: \mathcal{V} \to \mathcal{W}$ be a completely gauge contractive surjective map with kernel $\mathcal{J}$. Then $\pi$ is called a \textbf{complete gauge quotient map} if $\mathcal{W}$ is completely gauge isometric to $(\mathcal{V} / \mathcal{J}, q^{\mathcal{J}})$ via $\pi(x) \mapsto x + \mathcal{J}$, i.e. \[ \omega_n(\pi^{(n)}(x)) = q_n^{\mathcal{J}}(x + M_n(\mathcal{J})_h) = \inf \{ \nu_n(y) : y \in x + M_n(\mathcal{J})_h \} \]
for all $x \in \mathcal{W}$.
\end{definition}

\begin{theorem}
Let $(\mathcal{V}, \nu)$ be a matrix gauge $*$-vector space, and let $\widehat{V}$ be its unitization. Then $(\mathcal{V}, \nu)$ is completely gauge isometric to an operator system if and only if there exists a complete gauge quotient map $\pi: \widehat{V} \to V$. When this is the case, the unit of $\mathcal{V}$ is $\pi(0,1)$.
\end{theorem}

\begin{proof}
First, assume that there exists a complete gauge quotient map $\pi: \widehat{\mathcal{V}} \to \mathcal{V}$. Then, by Corollary \ref{cor: kernels of cp maps}, $\mathcal{V}$ is an operator system with unit $e = \pi(0,1)$ and $q = \nu^e$. Since $\pi$ is a complete gauge quotient map, $\nu = q = \nu^{e}$.

Now assume that $\mathcal{V}$ is an operator system with unit $e$. Let $\widehat{\mathcal{V}}$ be its unitization. By Proposition \ref{prop: unitize operator system}, the map $\pi(x,\lambda) = x + \lambda e$ is completely positive. We see that $\pi$ is surjective, since $x = \pi(x,0)$ for every $x \in \mathcal{V}$. It remains to check that $\pi$ is a complete gauge quotient map. To see this, first note that $\mathcal{J} := \ker(\pi) = \text{span} \{(-e,1)\}$. Thus
\begin{eqnarray}
q_n((x,0) + M_n(\mathcal{J})_h) & = & \inf\{ u_n((x + T \otimes e, -T)) : T=T^* \in M_n \} \nonumber \\
& \leq & u_n(x,0) = \nu_n(x). \nonumber
\end{eqnarray}
By the universal property of operator system quotients (see Proposition 3.6 of \cite{Kavruk}), the map $\tilde{\pi}: \widehat{\mathcal{V}} / \mathcal{J} \to \mathcal{V}$ given by $(x,\lambda) + \mathcal{J} \mapsto \pi(x,\lambda) = x + \lambda e$ is completely positive. Since $\tilde{\pi}$ is completely gauge contractive with respect to the order unit gauges on $\widehat{\mathcal{V}} / \mathcal{J}$ and $\mathcal{V}$, we have $\nu_n(x) \leq q_n((x,0) + M_n(\mathcal{J})_h)$. Thus $q = \nu$. So $\pi$ is a complete gauge quotient map.
\end{proof}


\section{Webster-Winkler duality for matrix gauge $*$-vector spaces}

In this section, we present a Webster-Winkler type duality theory for matrix gauge $*$-vector spaces. We show that the category of matrix gauge $*$-vector spaces is naturally dual to the category of pointed compact matrix convex sets, introduced in \cite{KennedyManorKim}. However, this duality is simpler than the one presented in \cite{KennedyManorKim}, which requires additional technical conditions which we will not need. These extra conditions are needed to develop a Webster-Winkler type duality theory for Werner's non-unital operator systems, which differ from matrix gauge $*$-vector spaces as we have already noted. We will summarize the differences between these two perspectives at the end of this section.

Many ideas in this section are the same as those summarized in Section 4.6 of \cite{FullerHartzLupiniRectangular18}, which relies on the results of an earlier preprint of the present paper. Thus, our main contribution here is filling in the details and relating the ideas to the recent paper \cite{KennedyManorKim}.

We first recall the definition of a matrix convex set, introduced originally by Webster and Winkler in \cite{webster1999krein}.

\begin{definition}
Let $\mathcal{V}$ be a complex vector space. We say that a sequence $\mathcal{J} = (J_n)$ where $J_n \subseteq M_n(\mathcal{V})$ for each $n$ is called a \textbf{ matrix convex set} if whenever $\alpha_i \in M_{n_i,n}$ and $x_i \in J_{n_i}$ for each $i \in \{1,2,\dots,k\}$ and $\sum \alpha_i^* \alpha_i = I_n$ we have $\sum \alpha_i^* x_i \alpha_i \in J_n$. If $\mathcal{V}$ is a locally convex topological vector space and each $J_n$ is compact, we call $\mathcal{J}$ a \textbf{compact} matrix convex set. 
\end{definition}

Expressions of the form $\sum \alpha_i^* x_i \alpha_i$ as described above are called \textbf{matrix convex combinations}. The natural morphisms between matrix convex sets are matrix affine maps, which we define now.

\begin{definition}
Let $\mathcal{J} = (J_n)$ and $\mathcal{K} = (K_n)$ be matrix convex sets. A sequence of functions $\{F_n: J_n \to K_n\}$, compactly denoted as $F: \mathcal{J} \to \mathcal{K}$, is called a \textbf{matrix affine map} provided that whenever $\alpha_i \in M_{n_i,n}$ and $x_i \in J_{n_i}$ for each $i \in \{1,2,\dots,k\}$ and $\sum \alpha_i^* \alpha_i = I_n$ we have
\[ F_n( \sum \alpha_i^* x_i \alpha_i ) = \sum \alpha_i^* F_{n_i}(x_i) \alpha_i. \]
If each $F_n$ is a homeomorphism and $F^{-1} = (F_n^{-1})$ is matrix affine, then $F$ is called a \textbf{matrix affine homeomorphism}.
\end{definition}

Webster and Winkler showed that the category of compact matrix convex sets is in duality with the category of operator systems. In the following, we let $A(\mathcal{J})$ denote the set of continuous matrix affine functions $\theta: \mathcal{J} \to \mathbb{C}$, where $\mathbb{C}$ is regarded as a matrix convex set whose $n\textsuperscript{th}$ level is given by $M_n$. We regard a matrix $\Theta=(\theta_{ij}) \in M_n(A(\mathcal{J}))$ as a matrix affine function $\Theta: \mathcal{J} \to M_n$ by setting $\Theta(x)_{ij} = (\theta_{ij}(x))$. With these identifications, an element $F \in M_n(A(\mathcal{J}))$ is considered to be positive if the matrix $\Theta_k(x)$ is positive semidefinite for all $k \in \mathbb{N}$ and $x \in J_k$. We let $e \in A(\mathcal{J})$ denote the matrix affine map $e_n(x)=I_n$ for all $x \in J_n$. With these identifications, $A(\mathcal{J})$ may be regarded as an operator system with unit $e$. See Section 3 of \cite{webster1999krein} for more details.

\begin{theorem}[Webster-Winkler, Proposition 3.5 of \cite{webster1999krein}] \label{thm: Webseter Winkler}
Let $\mathcal{S}$ be an operator system, and let $\mathcal{J}$ be compact matrix convex set. 
\begin{enumerate}
    \item The matrix state space $\text{UCP}(\mathcal{S}) = (\text{UCP}_n(\mathcal{S}))$, where $\text{UCP}_n(\mathcal{S}) := \{ \phi: \mathcal{S} \to M_n : \phi \text{ is ucp} \} $ for every $n \in \mathbb{N}$, is a compact matrix convex set with respect to the weak-$*$ topology.
    
    \item The matrix affine functions $A(\mathcal{J})$ is an operator system with unit $e \in A(\mathcal{J})$ and matrix ordering given by the positive matrix affine functions in $M_n(A(\mathcal{J}))$.
    
    \item For any operator system $\mathcal{S}$, the mapping $s \mapsto \hat{s}$, where $\hat{s} \in A(\text{UCP}(\mathcal{S}))$ is given by $\hat{s}(\phi)=\phi(s)$, is a unital complete order isomorphism from $\mathcal{S}$ to $A(\text{UCP}(\mathcal{S}))$.
    
    \item For any compact matrix convex set $\mathcal{J}$, the mapping $x \mapsto \hat{x}$, where $\hat{x} \in \text{UCP}_n(A(\mathcal{J}))$ is given by $\hat{x}(\Theta)=\Theta(x)$, is a matrix affine homeomorphism from $\mathcal{J}$ to $\text{UCP}(A(\mathcal{J}))$.
\end{enumerate}
\end{theorem}

We call the duality expressed in Theorem \ref{thm: Webseter Winkler} between operator systems and matrix convex sets \textbf{Webster-Winkler duality}. We now seek to develop a similar duality for matrix gauge $*$-vector spaces and pointed matrix convex sets. The following definition essentially comes from \cite{KennedyManorKim}, although we will require fewer restrictions in our definition. 

\begin{definition}
Let $V$ be a complex vector space and let $\mathcal{J}$ be a matrix convex set in $V$ and let $z \in J_1$. We call the pair $(\mathcal{J},z)$ a \textbf{pointed matrix convex set}. When $V$ is a locally convex topological vector space and $\mathcal{J}$ is compact, we call $(\mathcal{J},z)$ a \textbf{compact} pointed matrix convex set.
\end{definition}

In the Webster-Winkler type duality we will develop, pointed matrix convex sets correspond to the completely gauge contractive maps on a matrix gauge $*$-vector space. The next proposition shows that the set of completely gauge contractive maps is indeed a compact pointed matrix convex set.

\begin{proposition} \label{prop: CGC homeomorphic to UCP}
Let $(\mathcal{V}, \nu)$ be a matrix gauge $*$-vector space. For each $n \in \mathbb{N}$, let $\text{CGC}_n(\mathcal{V})$ denote the set of completely gauge contractive linear maps $\phi: \mathcal{V} \to M_n$. Then $\text{CGC}(\mathcal{V}) = (\text{CGC}_n(\mathcal{V}))$ is a compact matrix convex set, and the pair $(\text{CGC}(\mathcal{V}),0)$ is a pointed matrix convex set, where $0$ denotes the zero map. Moreover, if $\widehat{\mathcal{V}}$ denotes the unitization of $\mathcal{V}$, then the restriction map from $\text{UCP}(\widehat{\mathcal{V}})$ to $\text{CGC}(\mathcal{V})$ is a matrix affine homeomorphism.
\end{proposition}

\begin{proof}
Let $\varphi: \mathcal{V} \to M_n$ be an element of $\text{CGC}_n(\mathcal{V})$. By Corollary \ref{cor: Unique unital extension op sys} there exists a unital completely positive extension $\widetilde{\varphi}: \widehat{\mathcal{V}} \to M_n$. It follows that the restriction map $\pi: \text{UCP}_n(\widehat{\mathcal{V}}) \to \text{CGC}_n(\mathcal{V})$ given by $\pi(\psi) = \psi |_{\mathcal{V}}$ is surjective. Since $\text{CGC}_n(\mathcal{V})$ is the continuous image of the compact set $\text{UCP}_n(\widehat{\mathcal{V}})$, it is compact. Since every completely gauge contractive map has a unique unital extension to the unitization, the restriction map is bijective. It follows that the restriction map is a homeomorphism (c.f. Theorem 26.6 of \cite{MunkresTopology}).
\end{proof}

We now define an analog of the set of matrix affine maps $A(\mathcal{J})$ suitable for the setting of pointed matrix convex sets. Our definition is similar to Definition 3.10 of \cite{KennedyManorKim}, in that the vector space we define is identical to the one in \cite{KennedyManorKim}. However, we will equip this vector space with a matrix gauge $*$-vector space structure.

\begin{definition}
Let $(\mathcal{J},z)$ be a compact pointed matrix convex set in a locally convex topological vector space $V$. Let $A(\mathcal{J},z)$ denote the set of matrix affine functions $\theta \in A(\mathcal{J})$ which satisfy
\[ \theta_n(I_n \otimes z) = 0 \quad \text{for all } \quad n \in \mathbb{N}.\] For each $\Theta \in M_k(A(\mathcal{J},z))$, we define 
\[ \nu_k^{\mathcal{J}}(\Theta) := \inf \{ t > 0: \Theta_n(x) \leq t I_n \otimes I_k \text{ for all } n \in \mathbb{N} \text{ and } x \in J_n \} \]
and we let $\nu^{\mathcal{J}}$ denote the sequence $\{\nu_n^{\mathcal{J}}\}$.
\end{definition}

We now show that $\nu^{\mathcal{J}}$ is a matrix gauge on $A(\mathcal{J},z)$. In fact, $\nu^{\mathcal{J}}$ is just the restriction of the unique matrix gauge $\nu^e$ on $A(\mathcal{J})$ to the subspace $A(\mathcal{J},z)$.

\begin{proposition} \label{prop: A(J,z) subspace of A(J)}
Let $(\mathcal{J},z)$ be a compact pointed matrix convex set in a locally convex topological vector space $V$. Then $(A(\mathcal{J},z),\nu^{\mathcal{J}})$ is a matrix gauge $*$-vector space. Moreover, $\nu^{\mathcal{J}}$ is equal to the restriction of the matrix gauge $\nu^e$ on the operator system $A(\mathcal{J})$.
\end{proposition}

\begin{proof}
It is clear that $A(\mathcal{J},z)$ is a self-adjoint vector subspace of the operator system $A(\mathcal{J})$. If we can show that $\nu^{\mathcal{J}}$ is equal to the restriction of the matrix gauge $\nu^e$ on the operator system $A(\mathcal{J})$, it will follow that $\nu^{\mathcal{J}}$ is a matrix gauge on $A(\mathcal{J},z)$. If $\Theta \in M_n(A(\mathcal{J}))$, then
\[ \nu_n^e(\Theta) = \inf \{ t>0 : \Theta \leq t I_n \otimes e \}. \]
By the definition of the order structure on $A(\mathcal{J})$, we have $\Theta \leq tI_n \otimes e$ if and only if $\Theta_k(x) \leq t I_n \otimes e_k(x) = tI_n \otimes I_k$ for all $k \in \mathbb{N}$ and $x \in J_k$. Therefore $\nu_n^e(\Theta) = \nu_n^{\mathcal{J}}(\Theta)$ for every $\Theta \in M_n(A(\mathcal{J},z))$.
\end{proof}

We can now state our analog of Webster-Winkler duality for matrix gauge $*$-vector spaces. The results will follow easily from the identifications of $A(\mathcal{J},z)$ as a subspace of $A(\mathcal{J})$ and the matrix affine homeomorphism between $\text{UCP}(\widehat{\mathcal{V}})$ and $\text{CGC}(\mathcal{V})$.

\begin{theorem}
Let $(\mathcal{V}, \nu)$ be a matrix gauge $*$-vector space, and let $(\mathcal{J}, z)$ be a compact pointed matrix convex set. 
\begin{enumerate}
    \item The completely gauge contractive maps $\text{CGC}(\mathcal{V}) = (\text{CGC}_n(\mathcal{V}))$ together with the zero map form a pointed compact matrix convex set with respect to the weak-$*$ topology.
    
    \item The sequence $\nu^{\mathcal{J}}$ is a matrix gauge on the $*$-vector space $A(\mathcal{J},z)$.
    
    \item For any matrix gauge $*$-vector space $\mathcal{V}$, the mapping $v \mapsto \hat{v}$, where $\hat{v} \in A(\text{CGC}(\mathcal{V}),0)$ is given by $\hat{v}(\phi)=\phi(v)$, is a complete gauge isometry.
    
    \item For any compact matrix convex set $\mathcal{J}$, the mapping $x \mapsto \hat{x}$, where $\hat{x} \in \text{CGC}_n(A(\mathcal{J},z))$ given by $\hat{x}(\Theta)=\Theta(x)$, is a matrix affine homeomorphism from $\mathcal{J}$ to $\text{CGC}(A(\mathcal{J},z))$ sending $z$ to the zero map.
\end{enumerate}
\end{theorem}

\begin{proof}
Statement (1) follows from Proposition \ref{prop: CGC homeomorphic to UCP}, and statement (2) follows from Proposition \ref{prop: A(J,z) subspace of A(J)}.

We now consider statement (3). The mapping $v \mapsto \hat{v}$ is well-defined since $\hat{v}(0)=0$ for all $v \in \mathcal{V}$. Let $\mathcal{J} = \text{CGC}(\mathcal{V})$. We must show that the mapping $v \mapsto \hat{v}$ is onto and that $\nu_n(v) = \nu_n^{\mathcal{J}}(\hat{v})$ for each $v \in M_n(\mathcal{V})$. By Proposition \ref{prop: A(J,z) subspace of A(J)}, the inclusion map $i: A(\mathcal{J},0) \to A(\mathcal{J})$ is a complete gauge isometry identifying the unique matrix gauge $\nu^e$ on $A(\mathcal{J})$ with $\nu^{\mathcal{J}}$ on $A(\mathcal{J},0)$. By Proposition \ref{prop: CGC homeomorphic to UCP}, $\mathcal{J}$ is matrix affine homeomorphic to $\text{UCP}(\widehat{\mathcal{V}})$, yielding a unital complete order embedding $\pi: A(\mathcal{J}) \to A(\text{UCP}(\widehat{\mathcal{V}}))$. By Theorem \ref{thm: Webseter Winkler}, the mapping $j:s \mapsto \hat{s}$ is a unital complete order isometry of $\widehat{\mathcal{V}}$ onto $A(\widehat{\mathcal{V}})$, and hence a complete gauge isometry with respect to the order unit matrix gauges $\nu^e$. The result follows from checking that the gauge isometric inclusion of $\mathcal{V}$ into $\widehat{\mathcal{V}}$ coincides with the mapping $v \mapsto j^{-1} \circ \pi \circ i(\hat{v})$.

We conclude by considering statement (4). It is clear that $x \mapsto \hat{x}$ is well-defined and that $\hat{z}=0$. By Theorem \ref{thm: Webseter Winkler}, the mapping $x \mapsto \hat{x}$ is a matrix affine homeomorphism from $\mathcal{J}$ to $\text{UCP}(A(\mathcal{J}))$. By Proposition \ref{prop: CGC homeomorphic to UCP}, the restriction map induces a matrix affine homeomorphism from $\text{UCP}(\widehat{A(\mathcal{J},z)})$ to $\text{CGC}(A(\mathcal{J},z))$. The result will follow if $\widehat{A(\mathcal{J},z)} = A(\mathcal{J})$. From Proposition \ref{prop: A(J,z) subspace of A(J)}, we have that $A(\mathcal{J},z)$ is a gauge isometric subspace of $A(\mathcal{J})$. By Theorem \ref{thm: Unitization}, we have $\widehat{A(\mathcal{J},z)} \subseteq A(\mathcal{J})$ since the unitization is unique. Therefore we only need to check that the inclusion $A(\mathcal{J},z) \subseteq A(\mathcal{J})$ is codimension one. This follows from the identification of $A(\mathcal{J})$ with the continuous affine functions on $J_1$ (c.f. Section 3 of \cite{webster1999krein}) and the observation that if $\theta: J_1 \to \mathbb{C}$ is affine, then $\theta = \theta' + \theta(z)e$ where $e(x)=1$ for all $x \in J_1$ and $\theta' \in A(\mathcal{J},z)$.
\end{proof}

\begin{remark}
\emph{There are two significant differences between the Webster-Winkler duality presented here and the one presented in \cite{KennedyManorKim}. First, the authors of \cite{KennedyManorKim} define a \textit{pointed NC convex set} rather than a pointed matrix convex set. NC convex sets are families of convex spaces which include levels corresponding to infinite cardinalities, in addition to the finite cardinalities of a matrix convex set. We have not considered NC convex sets here, but we suspect that our Webster-Winkler duality could be defined in terms of NC convex sets if desired.}

\emph{The other significant difference between our Webster-Winkler duality and the one of \cite{KennedyManorKim} is that we do not require that every quasistate (i.e. completely positive completely contractive map) on the matrix ordered operator space $A(\mathcal{J},z)$ arise from evaluation at a point of $\mathcal{J}$ (see Definitions 3.12 and 3.14 of \cite{KennedyManorKim} for comparison). That assumption is valid if and only if $A(\mathcal{J},z)$ is a gauge maximal subspace of the operator system $A(\mathcal{J})$, since the completely gauge contractive maps coincide with the completely positive completely contractive maps in that case only (by Corollary \ref{cor: extensions of cpcc and gague maimal}).}
\end{remark}



\section{Duals of operator systems} \label{Sec: duals of operator systems}

We conclude with an application of matrix gauges. We will demonstrate how to represent duals and preduals of operator systems as concrete ordered operator spaces using matrix gauges. We begin with a definition.

\begin{definition}
Let $\mathcal{V}$ be a complex $*$-vector space, $\mathcal{C} = (C_n)$ be a matrix ordering on $\mathcal{V}$, and let $\hat{e}: \mathcal{V} \to \mathbb{C}$ be a stictly positive linear functional, i.e. $\hat{e}(x) > 0$ for all nonzero $x \geq 0$. Let $\mathcal{B} = (B_n)$, where 
\[ B_n := \{ x \in C_n : I_n \otimes \hat{e} (x) = I_n \}. \]
We call the triple $(\mathcal{V}, \mathcal{C}, \mathcal{B})$ a \textbf{noncommutative base space} if it satisfies the following conditions:
\begin{enumerate}
    \item For every $x \in M_n(\mathcal{V})_h$, there exists $t > 0$ and $b \in B_n$ such that $x \leq t b$.
    \item If $x \in M_n(\mathcal{V})_h$ and for every $t > 0$ there exists $b \in B_n$ such that $x \leq t b$, then $x \leq 0$.
\end{enumerate}
The sequence $\mathcal{B}$ is called the \textbf{matrix base} of $\mathcal{V}$.
\end{definition}

\begin{remark}
\emph{A noncommutative base space generalizes the notion of a \textit{base norm space}, a structure used to abstractly characterize the dual or predual of a function system. See Chapter 2 of \cite{AlfsenConvexBook71} for details.}
\end{remark}

\begin{remark}
\emph{It is easily verified that $\mathcal{B}$ is a matrix convex subset of $\mathcal{C}$. Our definition of $\mathcal{B}$ is in terms of the linear functional $\hat{e}$; however, it can be defined independently of this functional. Informally, $\mathcal{B}$ is a matrix convex subset of $\mathcal{C}$ whose first level is the intersection of an affine hyperplane with $C_1$ not containing zero. For brevity, we will not develop the details of this perspective here.}
\end{remark}

\begin{definition}
Let $(\mathcal{V}, \mathcal{C}, \mathcal{B})$ be a noncommutative base space. For each $n \in \mathbb{N}$, define $\nu_n^{\mathcal{B}}: M_n(\mathcal{V})_h \to [0,\infty)$ by
\[ \nu_n^{\mathcal{B}}(x) := \inf \{t > 0 : x \leq t b \text{ for some } b \in B_n \}. \]
We let $\nu^{\mathcal{B}}$ denote the sequence $(\nu_n^{\mathcal{B}})$.
\end{definition}

The following theorem shows that $\nu^{\mathcal{B}}$ is a matrix gauge, and that the base $\mathcal{B}$ can be specified from the gauge. This matrix gauge will allow us to construct operator representations of noncommutative base spaces.

\begin{theorem}
Let $(\mathcal{V}, \mathcal{C}, \mathcal{B})$ be a noncommutative base space. Then
\begin{enumerate}
    \item $\nu^{\mathcal{B}}$ is a proper matrix gauge on $\mathcal{V}$.
    \item $\mathcal{C} = \mathcal{C}^{\nu}$ where $\nu = \nu^{\mathcal{B}}$.
    \item We have $\nu_1^{\mathcal{B}}(x)=\hat{e}(x)$ for all $x \in C_1$.  The linear functional $\hat{e}$, and hence the sequence $\mathcal{B}$, is uniquely specified by these relations.
\end{enumerate}
Moreover, we have $\nu^{max} = \nu^{\mathcal{B}}$, where $\nu^{max}$ is the maximal matrix gauge inducing normal matrix ordered operator space $(\mathcal{V}, \mathcal{C}, \|\cdot\|^{\nu^{\mathcal{B}}})$.
\end{theorem}

\begin{proof}
We first prove (2), i.e. we show that $x \in C_n$ if and only if $\nu_n^{\mathcal{B}}(-x) = 0$. If $x \in C_n$, then $x + t b \geq 0$ for every $t > 0$ and $b \in B_n$, since $B_n \subseteq C_n$. Thus $\nu_n^{\mathcal{B}}(-x)=0$. On the other hand, suppose $x \in M_n(\mathcal{V})_h$ and $\nu_n^{\mathcal{B}}(-x) = 0$. Then for every $t > 0$ there exists $b \in B_n$ such that $-x \leq t b$. It follows from the definition of a noncommutive base space that $-x \leq 0$, and hence $x \in C_n$.

We now prove (1). Let $n \in \mathbb{N}$. We will first show that $\nu_n^{\mathcal{B}}$ is a proper gauge. Let $x \in M_n(\mathcal{V})_h$ and let $r > 0$. Then $rx \leq t b$ for some $b \in B_n$ if and only if $x \leq (t/r) b$ for some $b \in B_n$. We conclude that $\nu_n^{\mathcal{B}}(rx) = r \nu_n^{\mathcal{B}}(x)$. Next, let $x,y \in M_n(\mathcal{V})_h$ and suppose that $x \leq t b$ and $x \leq t' b'$ for some $b,b' \in B_n$. Then
\[ x + y \leq tb + t' b' = (t+t') \left( \frac{t}{t+t'} b + \frac{t'}{t+t'} b' \right). \]
Since $\frac{t}{t+t'} b + \frac{t'}{t+t'} b' \in B_n$ we conclude that $\nu_n^{\mathcal{B}}(x+y) \leq \nu_n^{\mathcal{B}}(x) + \nu_n^{\mathcal{B}}(y)$. The fact that $\nu_n^\mathcal{B}$ is proper is immediate from (2). Therefore each $\nu_n^{\mathcal{B}}$ is a proper gauge. We now check that $\nu^{\mathcal{B}}$ is a matrix gauge. Let $x \in M_n(\mathcal{V})_h$ and $y \in M_k(\mathcal{V})_h$. Suppose that $t > \nu_n^{\mathcal{B}}(x)$ and $r > \nu_k^{\mathcal{B}}(y)$. Then $x \leq t b$ for some $b \in B_n$ and $y \leq r b'$ for some $b' \in B_k$. It follows that $x \oplus y \leq \max(t,r) (b \oplus b')$. Since $b \oplus b' \in B_{n+k}$, we have $\nu_{n+k}^{\mathcal{B}}(x \oplus y) \leq \max(t,r)$. Now suppose that $s < \max(\nu_n^{\mathcal{B}}(x), \nu_k^{\mathcal{B}}(y))$. Without loss of generality, we may assume $s < \nu_n^{\mathcal{B}}(x) = \max(\nu_n^{\mathcal{B}}(x), \nu_k^{\mathcal{B}}(y))$. If $x \oplus y \leq s b''$ for some $b'' \in B_{n+k}$, then $x \leq s P b'' P^*$, where $P = [I_n 0_k]$ is the compression to the upper left $n \times n$ block. Since $P b'' P^* \in B_n$, we have $\nu_n^{\mathcal{B}}(x) \leq s$, a contradiction. We conclude that $\nu_{n+k}(x \oplus y) = \max(\nu_n^{\mathcal{B}}(x), \nu_k^{\mathcal{B}}(y))$. Finally, let $x \in M_n(\mathcal{V})_h$ and let $\alpha \in M_{n,k}$ with $\|\alpha^* \alpha\| \leq 1$. Suppose that $x \leq r b$ for some $b \in B_n$. Then
\[ \alpha^* x \alpha \leq r \alpha^* b \alpha \leq r \left( \alpha^* b \alpha + (I - \alpha^* \alpha)^{1/2} b (I - \alpha^* \alpha)^{1/2} \right). \]
Since $\alpha^* b \alpha + (I - \alpha^* \alpha)^{1/2} b (I - \alpha^* \alpha)^{1/2} \in B_k$, $\nu_k^{\mathcal{B}}(\alpha^* x \alpha) \leq r$. We conclude that $\nu_k^{\mathcal{B}}(\alpha^* x \alpha) \leq \nu_n^{\mathcal{B}}(x)$. A rescaling argument shows that $\nu_k^{\mathcal{B}}(\alpha^* x \alpha) \leq \|\alpha\|^2 \nu_n^{\mathcal{B}}(x)$ for general $\alpha \in M_{n,k}$. Therefore we conclude that $\nu^{\mathcal{B}}$ is a matrix gauge.

We now check statement (3). Suppose that $x \in C_1$. We claim that $\nu_1^{\mathcal{B}}(x) = \hat{e}(x)$. To see this, let $b = x/\hat{e}(x)$. Then $\hat{e}(b)=1$, so $b \in B_1$. Since $x = \hat{e}(x) b$, we have $\nu_1^{\mathcal{B}}(x) \leq \hat{e}(x)$. Now if $\nu_1^{\mathcal{B}}(x) < \hat{e}(x)$, then there exists $b' \in B_1$ and $r < \hat{e}(x)$ such that $x \leq r b'$. Since $\hat{e}$ is positive, $0 \leq r\hat{e}(b') - \hat{e}(x) = r - \hat{e}(x)$, implying that $\hat{e}(x) \leq r$, a contradiction. Thus $\nu_1^{\mathcal{B}}(x) = \hat{e}(x)$ for all $x \geq 0$. To see that this uniquely specifies $\hat{e}$, let $x \in \mathcal{V}_h$. Choose $t > 0$ and $b,b' \in B_1$ such that $x \leq tb$ and $-x \leq tb'$. Then $\pm x \leq t(b+b')$. Since
\[ x = \frac{1}{2}(x + t(b+b')) - \frac{1}{2}(-x + t(b+b')) \]
$x$ is a difference of positive elements. In particular,
\[ \hat{e}(x) = \frac{1}{2}\nu_1^{\mathcal{B}}(x + t(b+b')) - \frac{1}{2}\nu_1^{\mathcal{B}}(-x + t(b+b')). \]
So $\hat{e}(x)$ uniquely specified. The value of $\hat{e}$ on non-self adjoint elements is determined by considering real and imaginary parts.

For the final claim, recall that for each $x \in M_n(\mathcal{V})_h$,
\[ \nu_n^{max}(x) = \inf_{p \in C_n} \|x + p\|_{\nu^{\mathcal{B}}}. \]
By Theorem \ref{thm: max gauge}, $\nu_n^{\mathcal{B}}(x) \leq \nu_n^{max}(x)$, so we will show that $\nu_n^{max}(x) \leq \nu_n^{\mathcal{B}}(x)$. If $\nu_n^{\mathcal{B}}(x) \geq \nu_n^{\mathcal{B}}(-x)$, then
\[ \nu^{max}(x) \leq \|x+0\| = \max(\nu_n^{\mathcal{B}}(x), \nu_n^{\mathcal{B}}(-x)) = \nu_n^{\mathcal{B}}(x). \]
Now suppose that $\nu_n^{\mathcal{B}}(x) \leq \nu_n^{\mathcal{B}}(-x)$. For each $t > \nu_n^{\mathcal{B}}(x)$, there exists $b \in B_n$ such that $x \leq t b$. Let $p = tb - x$. Then
\[ \nu_n^{max}(x) \leq \|x + p\| = \max( \nu_n^{\mathcal{B}}( t b) + \nu_n^{\mathcal{B}}(-tb) ) = t. \]
It follows that $\nu_n^{max}(x) \leq \nu_n^{\mathcal{B}}(x)$. \end{proof}

\begin{remark}
Since $\nu^{\mathcal{B}} = \nu^{max}$, the set of completely positive completely contractive maps on an noncommutative base space coincides with the set of completely gauge contractive maps.
\end{remark}

\begin{definition} \label{defn: dual of base}
Let $(\mathcal{V}, \mathcal{C}, \mathcal{B})$ be a noncommutative base space. Let $\mathcal{V}^d$ denote the vector space of linear functionals bounded with respect to the norm $\|\cdot\|^{\nu^{\mathcal{B}}}$. Let $\mathcal{C}^d = (C_n^d)$, where
\[ C_n^d := \{ f \in M_n(\mathcal{V}^d) : f \text{ is completely positive} \} \]
where elements $f \in M_n(\mathcal{V}^d)$ are interpreted as $M_n$-valued maps in the usual way. We call the triple $(\mathcal{V}^d, \mathcal{C}^d, \hat{e})$ the \textbf{dual} of $(\mathcal{V}, \mathcal{C}, \mathcal{B})$.
\end{definition}

We now show that the dual of a noncommuative base space is an operator system.

\begin{theorem} \label{thm: dual of base is op sys}
Let $(\mathcal{V}, \mathcal{C}, \mathcal{B})$ be a noncommutative base space. Then its dual $(\mathcal{V}^d, \mathcal{C}^d, \hat{e})$ is an operator system.
\end{theorem}

\begin{proof}
We first check that $\hat{e} \in C_1^d$. Let $x \in M_n(\mathcal{V})_h$. Let $t > \nu_n^{\mathcal{B}}(x)$. Then there exists $b \in B_1$ such that $x \leq tb$. It follows that $I_n \otimes \hat{e}(x) \leq t I_n \otimes \hat{e}(b) = t I_n$. Since this holds for all $t > \nu_n^{\mathcal{B}}(x)$, we conclude that $\hat{e}(x) < \nu_1^{\mathcal{B}}(x) I_n$. It follows that $\hat{e}$ is completely gauge contractive, and hence $\hat{e}$ is completely positive and completely contractive.

Next we check that $\hat{e}$ is an order unit for $(\mathcal{V}^d, C_1^d)$. Let $f \in \mathcal{V}^d_h$, and assume that $f$ is contractive. Then for every $x \in C_1$, $|f(x)| \leq \|x\| = \nu(x) = \hat{e}(x)$. It follows that $(f + \hat{e})(x) \geq 0$. A rescaling argument shows that $f + \|f\|\hat{e} \in C_1^d$. So $\hat{e}$ is an order unit. Consequently, $\hat{e}$ is a matrix order unit.

We now check that $\hat{e}$ is an Archimedean. Let $f: \mathcal{V} \to M_n$ and suppose that $f + tI_n \otimes \hat{e}$ is completely positive for all $t > 0$. Let $x \in C_k$, fix $r > \nu_k^{\mathcal{B}}(x)$, and choose $b \in B_k$ such that $x \leq r b$. Since $\hat{e}$ is completely positive, $I_{nk} \otimes \hat{e}(x) \leq r I_{nk} \otimes \hat{e}(b) = r I_{nk}$. Hence, for all $t > 0$,
\[ 0 \leq f_k(x) + t I_{nk} \otimes \hat{e}(x) \leq f_k(x) + tr I_{nk}. \]
Since $I_{nk}$ is an Archimedean order unit in $M_{nk}$ and since $r > 0$ is fixed, $f_k(x) \geq 0$. This verifies that $\hat{e}$ is an Archimedean matrix order unit. We conclude that $(\mathcal{V}^d, \mathcal{C}^d, \hat{e})$ is an operator system. \end{proof}

The above results prove that every noncommutative base space, endowed with the operator space structure it inherits from the gauge $\nu^{\mathcal{B}}$, is the predual of an operator system. We now show that the operator space dual of an operator system $\mathcal{S}$ is a noncommutative base space $\mathcal{S}^d$. The matrix norm on $\mathcal{S}^d$ induced by its gauge is generally different from the operator space dual given by
\[ \|f\|_n^{d} = \|f\|_{cb} = \sup \{ \|f^{(k)}(x)\| : \|x\|_k \leq 1, k \in \mathbb{N} \}, \]
but we will see that these norms agree on the positive cone of $\mathcal{S}^d$. We will also see that the dual of $\mathcal{S}^d$, as defined in Definition \ref{defn: dual of base}, agrees with the operator space double dual of $\mathcal{S}$. We begin by defining a noncommutative base space structure on the operator space dual of an operator system.

\begin{definition}
Let $(\mathcal{S}, \mathcal{C}, e)$ be an operator system. Let $\mathcal{S}^d$ denote the vector space of linear functionals bounded with respect to the matrix norm $\|\cdot\|^{e}$. Let $\mathcal{C}^d = (C_n^d)$, where
\[ C_n^d := \{ f \in M_n(\mathcal{V}^d) : f \text{ is completely positive} \} \]
with elements $f \in M_n(\mathcal{V}^d)$ interpreted as $M_n$-valued maps in the usual way. Let $\mathcal{B} = (B_n)$ denote the matrix states
\[ B_n := \{f \in \mathcal{C}_n^d : f(e) = I_n \}. \]
We call the triple $(\mathcal{S}^d, \mathcal{C}^d, \mathcal{B})$ the \textbf{dual} of $(\mathcal{S}, \mathcal{C}, e)$.
\end{definition}


\begin{theorem} \label{thm: opSys Duality}
Let $(\mathcal{S}, \mathcal{C}, e)$ be an operator system. Then
\begin{enumerate}
    \item The dual $(\mathcal{S}^d, \mathcal{C}^d, \mathcal{B})$ is a noncommutative base norm space with $\hat{e}$ given by evaluation at $e$.
    \item The matrix norm $\|\cdot\|^{\nu^{\mathcal{B}}}$ agrees with the operator space dual norm $\|\cdot\|^d$ on elements of the matrix ordering $\mathcal{C}^d$.
    \item The operator system $(\mathcal{S}^{dd}, \mathcal{C}^{dd}, e)$, regarded as the dual of $(\mathcal{S}^d, \mathcal{C}^d, \mathcal{B})$, coincides with the operator space double dual of $(\mathcal{S}, \mathcal{C}, e)$, so that the natural embedding $\mathcal{S} \hookrightarrow \mathcal{S}^{dd}$ is a unital complete order embedding.
\end{enumerate}
\end{theorem}

\begin{proof}
We first prove statement (1). It is known that $\mathcal{C}^d$ is a matrix ordering on $\mathcal{S}^d$. Let $\hat{e}: \mathcal{S}^d \to \mathcal{C}$ be given by evaluation at $e$, i.e. $\hat{e}(f) = f(e)$. For any $f \in \mathcal{S}^d$, $\|f\| = |f(e)|$ (c.f. Proposition 3.6 of \cite{paulsen2002completely}). Thus, if $f \geq 0$ and $f \neq 0$, then $\hat{e}(f) > 0$. It follows from the definition of $\mathcal{B}$ that $(\mathcal{S}^d, \mathcal{C}^d, \mathcal{B})$ is a noncommutative base space.

We now prove statement (2). Suppose that $f \in C_n^d$. Then 
\[ \|f\|^{\nu^{\mathcal{B}}}_n = \nu^{\mathcal{B}}_n(f) = \inf \{ t > 0: f \leq t b \text{ for some } b \in B_n\}. \]
Since $f$ is completely positive, $\|f\|_{cb} = \|f(e)\|$ (using Proposition 3.6 of \cite{paulsen2002completely} again). Let $\varphi: \mathcal{S} \to \mathbb{C}$ be an element of $C_1^d$, i.e. a state on $\mathcal{S}$. Let $P = \|f(e)\|I_n - f(e)$. Then $P \geq 0$, and hence the map $g(x) = \varphi(x) P$ is an element of $C_n^d$. Define 
\[ b(x) := \frac{1}{\|f(e)\|} \left( f(x) + g(x) \right). \]
Then $b \in B_n$ and
\[ f \leq f + g = \|f(e)\| b. \]
So we have $\nu^{\mathcal{B}}(x) \leq \|f(e)\| = \|f\|_{cb} = \|f\|_n^d$. On the other hand, suppose $f \leq r b$ for some $r < \|f(e)\|$ and $b \in B_n$. Then
\[ 0 \leq f(e) \leq r b(e) = r I_n, \]
a contradiction, since this implies $\|f(e)\| \leq r$. We conclude that $\|f\|_n^d = \|f\|_n^{\nu^{\mathcal{B}}}$, proving the statement.

We now prove statement (3). It suffices to prove that the operator space dual norm on $\mathcal{S}^d$ is equivalent to the norm induced by the gauge, so that continuity with respect to one norm implies continuity with respect to the other. To this end, suppose that $f: \mathcal{S} \to \mathbb{C}$ is self-adjoint. Then $f = f_+ - f_-$, where $f_+, f_- \geq 0$ and $\|f\|_1^d = \|f_+\|_1^d + \|f_-\|_1^d$ (see Theorem 4 of \cite{EllisDuality64}). By (2), $\|f_+\|_1^d = \nu_1^{\mathcal{B}}(f_+)$ and $\|f_-\|_1^d = \nu_1^{\mathcal{B}}(f_-)$. Thus $\|f\|_1^d = \nu_1^{\mathcal{B}}(f_+) + \nu_1^{\mathcal{B}}(f_-)$. We claim that $\nu_1^{\mathcal{B}}(f) = \nu_1^{\mathcal{B}}(f_+)$ and $\nu_1^{\mathcal{B}}(-f) = \nu_1^{\mathcal{B}}(f_-)$. To see this, suppose that $f \leq tb$ for some $b \in \mathcal{B}$. Then $f = tb - (tb-f)$, and hence 
\[ \|f_+\|^d_1 + \|f_-\|^d_1 = \|f\|^d_1 \leq \|tb\|^d_1 + \|tb-f\|^d_1 = tb(e) + tb(e) - f(e) = 2t - \|f_+\|^d_1 + \|f_-\|^d_1. \] 
Thus $\|f_+\|^d_1 \leq t$, and we conclude that $\|f_+\|^d_1 \leq \nu_1^{\mathcal{B}}(f)$. Since $f_+/\|f_+\|^d_1 \in \mathcal{B}$, we also have $\nu_1^{\mathcal{B}}(f) \leq \|f_+\|^d_1$. Therefore $\nu_1^{\mathcal{B}}(f) = \|f_+\|^d_1$. A similar argument shows $\nu_1^{\mathcal{B}}(-f) = \|f_-\|^d_1$. Now
\[ \max( \nu_1^{\mathcal{B}}(f_+), \nu_1^{\mathcal{B}}(f_-)) \leq \nu_1^{\mathcal{B}}(f_+) + \nu_1^{\mathcal{B}}(f_-) \leq 2 \max(\nu_1^{\mathcal{B}}(f_+), \nu_1^{\mathcal{B}}(f_-)). \]
We conclude that $\|f\|_1^{\nu} \leq \|f\|_1^d \leq 2 \|f\|_1^{\nu}$. Hence $\|\cdot\|_1^d$ is equivalent to $\|\cdot\|_1^{\nu}$. So $\mathcal{S}^{dd}$ coincides with the standard double dual of $\mathcal{S}$. \end{proof}

\begin{corollary}
Let $(\mathcal{S}, \mathcal{C}, e)$ be an operator system with dual $(\mathcal{S}^d, \mathcal{C}^d, \mathcal{B})$. Then there exists a complete order embedding $\pi: \mathcal{S}^d \to B(H)$ with the property that for each $f \in C_n^d$, $\|\pi^{(n)}(f)\| = \|f\|_n^d$.
\end{corollary}

\begin{proof}
Take $\pi$ to be a complete gauge isometry with respect to $\nu^{\mathcal{B}}$. Then the result follows from Theorem \ref{thm: opSys Duality} and Corollary \ref{cor: Representations of matrix gauge spaces}.
\end{proof}

\bibliographystyle{plain}
\bibliography{references}

\begin{thebibliography}{10}

\bibitem{AlfsenConvexBook71}
Erik~M. Alfsen.
\newblock {\em Compact convex sets and boundary integrals}.
\newblock Ergebnisse der Mathematik und ihrer Grenzgebiete ; Bd. 57.
  Springer-Verlag, Berlin, 1971.

\bibitem{Blecher2007}
David~P. Blecher.
\newblock {\em Positivity in Operator Algebras and Operator Spaces}, pages
  27--71.
\newblock Birkh{\"a}user Basel, Basel, 2007.

\bibitem{BlecherKirkpartickNealWerner2007}
David~P. Blecher, Kay Kirkpartick, Matthew Neal, and Wend Werner.
\newblock Ordered involutive operator spaces.
\newblock {\em Positivity}, 11:497--510, 2007.

\bibitem{BRS90}
David~P Blecher, Zhong-Jin Ruan, and Allan~M Sinclair.
\newblock A characterization of operator algebras.
\newblock {\em Journal of Functional Analysis}, 89(1):188--201, 1990.

\bibitem{CHOIEffros1977}
Man-Duen Choi and Edward~G. Effros.
\newblock Injectivity and operator spaces.
\newblock {\em Journal of Functional Analysis}, 24(2):156--209, 1977.

\bibitem{ConnesSuijlekomI}
Alain Connes and Walter~D. van Suijlekom.
\newblock Spectral truncations in noncommutative geometry and operator systems.
\newblock {\em Communications in Mathematical Physics}, 383:2021–--2067,
  2021.

\bibitem{ConnesSuijlekomII}
Alain Connes and Walter~D. van Suijlekom.
\newblock Tolerance relations and operator systems.
\newblock {\em arXiv:2111.02903}, 2021.

\bibitem{DavidsonCstarBook}
Kenneth~R. Davidson.
\newblock {\em {$C^*$}-algebras by example}, volume~6 of {\em Fields Institute
  Monographs}.
\newblock American Mathematical Society, Providence, RI, 1996.

\bibitem{EFFROSWinkler1997}
Edward~G. Effros and Soren Winkler.
\newblock Matrix convexity: Operator analogues of the {B}ipolar and
  {H}ahn–{B}anach {T}heorems.
\newblock {\em Journal of Functional Analysis}, 144(1):117--152, 1997.

\bibitem{EllisDuality64}
A.~J. Ellis.
\newblock {The Duality of Partially Ordered Normed Linear Spaces}.
\newblock {\em Journal of the London Mathematical Society}, s1-39(1):730--744,
  01 1964.

\bibitem{FullerHartzLupiniRectangular18}
Adam~H Fuller, Michael Hartz, and Martino Lupini.
\newblock Boundary representations of operator spaces, and compact rectangular
  matrix convex sets.
\newblock {\em Journal of Operator Theory}, 79(1):139--172, 2018.

\bibitem{karn_2011}
Anil~K. Karn.
\newblock Order embedding of a matrix ordered space.
\newblock {\em Bulletin of the Australian Mathematical Society}, 84(1):10–18,
  2011.

\bibitem{Kavruk}
Ali~S Kavruk, Vern~I Paulsen, Ivan~G Todorov, and Mark Tomforde.
\newblock Quotients, exactness, and nuclearity in the operator system category.
\newblock {\em Advances in mathematics}, 235:321--360, 2013.

\bibitem{KennedyManorKim}
Matthew Kennedy, Se-Jin Kim, and Nicholas Manor.
\newblock {Nonunital Operator Systems and Noncommutative Convexity}.
\newblock {\em International Mathematics Research Notices}, 01 2022.
\newblock rnab349.

\bibitem{MunkresTopology}
James~R. Munkres.
\newblock {\em Topology}.
\newblock Prentice Hall, Inc., 2 edition, 2000.

\bibitem{Ng}
Chi-Keung Ng.
\newblock Operator subspaces of {$L(H)$} with induced matrix orderings.
\newblock {\em Indiana University Mathematics Journal}, 60(2):577--610, 2011.

\bibitem{Ozawa2008CAlgebrasAF}
Nathanial P. Brown~Narutaka Ozawa.
\newblock C*-algebras and finite-dimensional approximations.
\newblock 2008.

\bibitem{PaulsenTomforde}
Vern Paulsen and Mark Tomforde.
\newblock Vector spaces with an order unit.
\newblock {\em Indiana University Mathematics Journal}, pages 1319--1359, 2009.

\bibitem{paulsen2002completely}
Vern~I Paulsen.
\newblock {\em Completely bounded maps and operator algebras}, volume~78.
\newblock Cambridge University Press, 2002.

\bibitem{PaulsenTodorovTomforde2011OSS}
Vern~I. Paulsen, Ivan~G. Todorov, and Mark Tomforde.
\newblock Operator system structures on ordered spaces.
\newblock {\em Proceedings of the London Mathematical Society}, 102(1):25--49,
  2011.

\bibitem{RUAN1988}
Zhong-Jin Ruan.
\newblock Subspaces of {C$^*$}-algebras.
\newblock {\em Journal of Functional Analysis}, 76(1):217--230, 1988.

\bibitem{RUSSELL2017}
Travis~B. Russell.
\newblock Characterizations of ordered operator spaces.
\newblock {\em Journal of Mathematical Analysis and Applications},
  452(1):91--108, 2017.

\bibitem{webster1999krein}
Corran Webster and Soren Winkler.
\newblock The {K}rein-{M}ilman theorem in operator convexity.
\newblock {\em Transactions of the American Mathematical Society},
  351(1):307--322, 1999.

\bibitem{WERNER2002}
Wend Werner.
\newblock Subspaces of {$L(H)$} that are $*$-invariant.
\newblock {\em Journal of Functional Analysis}, 193(2):207--223, 2002.

\end{thebibliography}

\end{document}